\documentclass[12pt]{article}
\usepackage[paperwidth=21cm,paperheight=26cm, textwidth=16cm,textheight=20.6cm]{geometry}
\usepackage{amsmath,amsthm,amssymb,xspace,xcolor,hyperref}\hypersetup{colorlinks=true, linkcolor=black,citecolor=black,urlcolor=black}\usepackage{graphicx}
\numberwithin{equation}{section}

\theoremstyle{plain}
\newtheorem{theorem}{Theorem}

\newtheorem{claim}{Claim}[section]
\newtheorem{proposition}{Proposition}
\newtheorem{lemma}{Lemma}
\newtheorem{corollary}{Corollary}
\newtheorem{definition}{Definition}
\newtheorem{question}{Question}

\newcommand{\APC}{\ensuremath{{\sf APC_1}}\xspace}

\newcommand{\SB}{\ensuremath{{\sf S^1_2}}\xspace}

\newcommand{\PV}{\ensuremath{{\sf PV_1}}\xspace}
\newcommand{\pv}{\ensuremath{{\sf PV}}\xspace}

\newcommand{\EF}{\ensuremath{{\sf EF}}\xspace}
\newcommand{\WF}{\ensuremath{{\sf WF}}\xspace}

\newcommand{\Ptime}{\ensuremath{{\sf P}}\xspace}
\newcommand{\Ppoly}{\ensuremath{{\sf P/poly}}\xspace}

\newcommand{\NPtime}{\ensuremath{{\sf NP}}\xspace}
\newcommand{\NP}{\ensuremath{{\sf NP}}\xspace}
\newcommand{\coNP}{\ensuremath{{\sf coNP}}\xspace}

\newcommand{\Circuit}{\ensuremath{{\sf Circuit}}\xspace}

\newcommand{\SAT}{\ensuremath{{\sf SAT}}\xspace}

\newcommand{\MCSP}{\ensuremath{{\sf MCSP}}\xspace}
\newcommand{\GCSP}{\ensuremath{{\sf GCSP}}\xspace}

\newcommand{\LB}{\ensuremath{{\sf LB}}\xspace}
\newcommand{\LBtt}{\ensuremath{{\sf LB_{tt}}}\xspace}
\newcommand{\ttable}{\ensuremath{{\sf tt}}\xspace}

\begin{document}

\title{{\bf Towards P $\ne$ NP from Extended\\ Frege lower bounds}\footnote{A version of the paper with a less technical introduction appeared on ECCC.}} 
\author{J\'an Pich \\ \small University of Oxford \and Rahul Santhanam \\ \small University of Oxford}
\date{\small September 2023}
\maketitle

\begin{abstract} We prove that if conditions I-II (below) hold and 
there is a sequence of Boolean functions $f_n$ hard to approximate by p-size circuits such that p-size circuit lower bounds for $f_n$ do not have p-size proofs in Extended Frege system \EF, then $\Ptime\ne \NP$.

\begin{itemize}
\item[I.] $\SB$ proves that a concrete function in ${\sf E}$ is hard to approximate by subexponential-size circuits.
\item[II.] {\bf Learning from $\neg\exists$ OWF.} \SB proves that a p-time reduction transforms circuits breaking one-way functions to p-size circuits learning p-size circuits over the uniform distribution, with membership queries.
\end{itemize}
Here, \SB is Buss's theory of bounded arithmetic formalizing p-time reasoning.

Further, we show that any of the following assumptions implies that $\Ptime\ne\NP$, if \EF is not p-bounded: 
\begin{itemize}
\item[1.] {\bf Feasible anticheckers.} \SB proves that a p-time function generates anticheckers for \SAT.
\item[2.] {\bf Witnessing $\NP\not\subseteq\Ppoly$.} \SB proves that a p-time function witnesses an error of each p-size circuit which fails to solve \SAT.

\item[3.] {\bf OWF from $\NP\not\subseteq\Ppoly$ \& hardness of {\sf E}.} Condition I holds and \SB proves that a p-time reduction transforms circuits breaking one-way functions to p-size circuits computing \SAT.
\end{itemize}
The results generalize to stronger theories and proof systems. 
\end{abstract}


\section{Introduction}
The {\em proof complexity} approach to the \Ptime versus \NP problem, sometimes refereed to as the Cook-Reckhow program \cite{BP,Kpc}, proceeds by proving lower bounds on lengths of proofs of tautologies in increasingly powerful proof systems - $\NP\ne \coNP$ (and hence $\Ptime\ne\NP$) iff all propositional proof systems have hard sequences of tautologies that require superpolynomial proof size. A problem with the approach is that we do not know if we ever reach the point of proving a superpolynomial lower bound for all proof systems, if we focus only on concrete ones.

There are several results which mitigate this situation. For example, `lifting' theorems provide a method for deriving monotone circuit lower bounds from lower bounds for weak proof systems \cite{Ls}. In the algebraic setting, superpolynomial lower bounds for CNF tautologies in the Ideal Proof System imply that the Permanent does not have p-size arithmetic circuits \cite{GP}. Proving that a proof system $P$ is not p-bounded implies also that $\SAT\notin \Ppoly$ is consistent with a theory $T_P$ corresponding to the proof system $P$ via the standard correspondence in bounded arithmetic \cite{Kpc}. 

However, no generic way of deriving explicit circuit lower bounds for unrestricted Boolean circuits from proof complexity lower bounds for concrete propositional proof systems of interest has been discovered.\footnote{It is essentially a trivial observation that there is some proof system $P$ such that if $P$ is not p-bounded then $\NPtime\ne \coNP$. To see that, consider two cases: I. If $\NP\ne \coNP$, we can take arbitrary $P$; II. If $\NP=\coNP$, we can take a p-bounded proof system for $P$. 
It would be dramatically different to obtain such a connection for a concrete proof system.} In the present paper we connect this problem to several classical questions in complexity theory.

\subsection{Our contribution}

\subsubsection{Core idea: Self-provability of circuit upper bounds}

The initial point for our considerations is an observation about a potential `self-provability' of circuit upper bounds: As it turns out, under certain assumptions, the mere truth of $\NP\subseteq\Ppoly$ implies its own provability. To explain the self-provability we first consider a simple-to-state question about a witnessing of $\NP\not\subseteq\Ppoly$.
\medskip

\noindent {\bf Witnessing $\NP\not\subseteq\Ppoly$.} Fix a constant $k\ge 1$. Suppose that for each sufficiently big $n$, no circuit with $n$ inputs and size $n^k$ finds a satisfying assignment for each satisfiable propositional formula of size $n$, i.e. no $n^k$-size circuit solves the search version of \SAT. Can we witness this assumption `feasibly' by a p-time function $f$ such that for each sufficiently big $n$, for each $n^k$-size circuit $C$ with $n$ inputs and $\le n$ outputs, $f(C)$ outputs a formula $\phi$ of size $n$ together with a satisfying assignment $a$ of $\phi$ such that $\neg \phi(C(\phi))$?

Similar kinds of witnessing have been considered before in the literature, using diagonalization techniques \cite{GST,Abb,Hardins}. Indeed, Bogdanov, Talwar and Wan \cite{Hardins} call a similar feasible witnessing in the uniform setting a ``dreambreaker" (following Adam Smith) and show that such a feasible witnessing can be constructed. However, in our applications, it is crucial that the witnessing function finds an error on the input length of the given circuit assuming just that the circuit errs on the given input length, and we do not yet know an unconditional answer to the question in the previous paragraph.

A witnessing function $f$ from the penultimate paragraph exists under the assumption of the existence of a one-way function and a function in {\sf E} hard for subexponential-size circuits \cite{Pclba,MP}, see also Lemma \ref{l:owpwit}. Is it, however, possible to construct it without assuming more than the assumption we want to witness? Formally, we are asking if there is a p-time function $f$ such that for each big enough $n$ propositional formula $w^k_n(f)$, defined by \begin{equation*}\label{e:witness}w^k_n(f):=[\SAT_n(x,y)\rightarrow \SAT_n(x,C(x))]\vee [\SAT_n(f_1(C),f_2(C))\wedge\neg\SAT_n(f_1(C),C(f_1(C)))],\end{equation*} is a tautology. Here, $\SAT_n(x,y)$ is a p-time predicate saying that $x$ is an $n$-bit string encoding a propositional formula satisfied by assignment $y$, $C(z)$ says that free variables $C$ represent a circuit with $n$ inputs, $\le n$ outputs and size $n^k$ which outputs $C(z)$ on $z$, and $f(C)$ outputs a pair of strings $\left<f_1(C), f_2(C)\right>$. We do not specify the precise encoding of $w^k_n(f)$. The proof systems we work with simulate Extended Frege system \EF and are therefore strong enough to reason efficiently with any natural encoding of $w^k_n(f)$.
\medskip

If we had a function $f$ such that for some $n_0$ and all $n>n_0$, $w^k_n(f)$ would be a tautology, we could define an extension of \EF, denoted $\EF+w^k(f)$, such that $\EF+w^k(f)$ proofs are \EF-proofs which are, in addition, allowed to derive substitutional instances of $w^k_n(f)$, for $n>n_0$. 
\medskip

\noindent {\bf Self-provability of $\NP\subseteq\Ppoly$.} Assuming we have a p-time function $f$ such that $w^k_n(f), n>n_0,$ is a tautology, we can now turn to the self-provability of $\SAT\in\Ppoly$.

Suppose that $\forall n>1, \SAT_n\in\Circuit[n^{k'}]$, where $\Circuit[s(n)]$ stands for the set of all single-output circuits with $n$ inputs and size $\le s(n)$. Then, there is a sequence of circuits $C$ with $n$ inputs, $\le n$ outputs and size $n^k$ falsifying the right disjunct in $w^k_n(f)$, for some $k>k'$ and all $n>1$. Therefore, $\EF+w^k(f)$ admits p-size proofs of $\SAT_n(x,y)\rightarrow \SAT_n(x,C(x))$. That is, the mere validity of $\SAT_n\in\Circuit[n^{k'}]$ implies an efficient propositional provability of $\SAT_n\in\Ppoly$. The efficient provability of $\SAT_n(x,y)\rightarrow \SAT_n(x,C(x))$ further implies that $\EF+w^k(f)$ is p-bounded: To prove a tautology $\phi$ of size $n$ in $\EF+w^k(f)$ it suffices to check out that $\neg\SAT_n(\neg\phi,C(\neg\phi))$ (which implies that $\SAT_n(\phi,y)$ and $\phi$ hold). \medskip

\noindent {\bf First-order setting.} We can formulate also a first-order version of the observation from the previous paragraph. Consider Buss's theory of bounded arithmetic \SB formalizing p-time reasoning, cf. \S\ref{sba}. Let $W^k_{n_0}(f)$ denote a natural $\forall \Pi^b_1$-formalization of the statement ``$\forall n>n_0, w^k_n(f)$", see \S\ref{sba} for the definition of $\Pi^b_1$. By the correspondence between \SB and \EF, cf. \cite{Kpc}, if $W^k_{n_0}(f)$ was provable in \SB, for some p-time $f$, then tautologies $w^k_n(f)$, for $n\ge n_0$, would have p-size proofs in \EF.

\medskip
To summarize, we proved the following.

\begin{theorem}[Circuit complexity from proof complexity \& witnessing of $\NP\not\subseteq\Ppoly$]\label{t:core}\ \newline Let $k\ge 1$ be a constant.
\begin{itemize}
\item[1.] Suppose that there is a p-time function $f$ such that for each big enough $n$, $w^k_n(f)$ is a tautology. If $\EF+w^k(f)$ is not p-bounded, then $\SAT_n\notin\Circuit[n^{\epsilon k}]$ for infinitely many $n$.
\item[2.] Suppose that there is a p-time function $f$ such that for some $n_0$, $\SB\vdash W^k_{n_0}(f)$. If $\EF$ is not p-bounded, then $\SAT_n\notin\Circuit[n^{\epsilon k}]$ for infinitely many $n$.
\end{itemize}
In Items 1 and 2, $\epsilon>0$ is a universal constant (independent of $k$).
\end{theorem}

Notably, if for all $k\ge 1$ there is a p-time function $f^k$ such that for each big enough $n$, $w^k_n(f^k)$ is a tautology, then $\NP\not\subseteq {\sf SIZE}[n^{\ell}]$, for all $\ell\in\mathbb{N}$. This is because by Theorem \ref{t:core}, the existence of such functions $f^k$ and the assumption that for each $k$, $\EF+w^k(f^k)$ is not p-bounded imply $\NP\not\subseteq\Ppoly$. On the other hand, if for some $k$, $\EF+w^k(f^k)$ is p-bounded, then $\NP=\mathsf{coNP}$ and 
$\NP\not\subseteq {\sf SIZE}[n^{\ell}]$, for each $\ell\in\mathbb{N}$, by Kannan's lower bound \cite{Klb}. 

\begin{corollary}[Circuit lower bounds from witnessing]\label{c:clbfromwit}\ \newline If for all $k\ge 1$ there is a p-time function $f^k$ such that for each big enough $n$, $w^k_n(f^k)$ is a tautology, then $\NP\not\subseteq {\sf SIZE}[n^{\ell}]$, for all $\ell\in\mathbb{N}$.
\end{corollary}

The last observation tells us also that obtaining an unconditional reduction of non-uniform circuit lower bounds to proof complexity of concrete systems will be challenging.

\begin{proposition}\label{p:hardideal} If there is a proof system $R$ such that showing that $R$ is not p-bounded implies $\NP\not\subseteq\Ppoly$, then $\NP\not\subseteq {\sf SIZE}[n^{\ell}]$, for all $\ell\in\mathbb{N}$.
\end{proposition}

\noindent {\bf Restricting nonuniformity.} In Theorem \ref{t:core}, we can restrict the number of nonuniform bits in the concluded lower bounds by adapting formulas $w^k_n(f)$: Assume that the circuit $C$ includes a hardwired description of a fixed universal Turing machine $U$. Moreover, interpret $C$ as encoding an algorithm $A$ described by $\le \log n$ bits with $u(n)\le n^k$ nonuniform bits of advice $a(n)$. 
The algorithm $A$ and its nonuniform advice are described by free variables. We assume that $u$ is p-time. On each input $z\in \{0,1\}^n$, $C(z)$ uses $U$ to simulate the computation of $A$ on $z$ with access to $a(n)$ up to $n^k$ steps. That is, now the size of $C$ is $poly(n^k)$. Denote the resulting formulas by $w^{k,u}_n(f)$. If we have a p-time function $f$ which witnesses errors of $n^k$-time algorithms described by $\log n$ bits with $u(n)$ bits of advice attempting to solve the search version of $\SAT$, i.e. such that formulas $w^{k,u}_n(f)$ are tautologies for big enough $n$, we can define proof system $\EF+w^{k,u}(f)$. Further, we can define $\forall\Pi^b_1$ formulas $W^{k,u}_n(f)$ expressing $``\forall n>n_0,\ w^{k,u}_n(f)$." Denote by ${\sf Time}[n^k]/u(n)$ the class of problems solvable by uniform algorithms with $\le u(n)$ bits of nonuniform advice running in time $O(n^k)$. The proof of Theorem \ref{t:core} works in this case as well. 

\begin{corollary}[Circuit complexity from proof complexity \& witnessing of $\Ptime\ne\NP$]\label{c:core}\ \newline Let $k\ge 1$ be a constant and $u$ a p-time function such that $u(n)\le n^k$.
\begin{itemize}
\item[1.] Suppose that there is a p-time function $f$ such that for each big enough $n$, $w^{k,u}_n(f)$ is a tautology. If $\EF+w^{k,u}(f)$ is not p-bounded, then $\SAT\notin {\sf Time}[n^{\Omega(k)}]/u(n)$.
\item[2.] Suppose that there is a p-time function $f$ such that for some $n_0$, $\SB\vdash W^{k,u}_{n_0}(f)$. If $\EF$ is not p-bounded, then $\SAT\notin{\sf Time}[n^{\Omega(k)}]/u(n)$.
\end{itemize}
\end{corollary}

The significance of Corollary \ref{c:core} is that in the uniform setting, a similar kind of feasible witnessing is known to exist using diagonalization techniques \cite{GST,Hardins}. It is unclear whether diagonalization techniques will suffice to establish that $w^{k,u}_n(f)$ is a tautology for large enough $n$ in some concrete proof system, but there is at least a strong motivation for considering the question, given its implications for deriving strong computational complexity lower bounds from proof complexity lower bounds.
\medskip

\def\provfromtruth{
Note that, similarly, it can be derived that if there is a function $f(x)$ outputting a triple $\left<f_1(x),f_2(x),f_3(x)\right>$ such that {\sf ZFC} proves that ``for each algorithm $A$, \begin{equation*} \Ptime=\NP \vee [\SAT_{f_3(A)}(f_1(A),f_2(A))\wedge\neg\SAT_{f_3(A)}(f_1(A),A(f_1(A)))]\vee R(A) ,"\end{equation*} where $R(A)$ says that $A(f_1(A))$ runs in time $>f_3(A)^k$, then the mere truth of $\SAT\in {\sf Time}[n^k]$ implies {\sf ZFC}-provability of $\Ptime=\NP$.}
\medskip

\noindent {\bf Nonuniform witnessing.} If we allow $f$ to be nonuniform, we obtain a version of formulas $w_n^k(f)$ which are unconditionally tautological. This follows from a theorem of Lipton and Young \cite{LY}, who showed that for each sufficiently big $n$ and each Boolean function $f$ with $n$ inputs which is hard for circuits of size $s^{3}$, $s\ge n^3$, there is a set $A^{f,s}_n\subseteq \{0,1\}^n$ of size $poly(s)$ such that no $s$-size circuit computes $f$ on $A^{f,s}_n$. The set $A^{f,s}_n$ is the set of {\em anticheckers} of $f$ w.r.t. $s$. Let $n$ be sufficiently big and $\alpha^s_n$ be tautologies defined by
\begin{equation*}\label{e:witness}\alpha^s_n:=\big(\SAT_n(x,y)\rightarrow \SAT_n(x,B(x))\big)\vee \big(\bigvee_{z\in A}C(z)\ne\SAT_n(z)\big),\end{equation*}
where $\SAT_n(z)\in \{0,1\}$ is such that $\SAT_n(z)=1\Leftrightarrow \exists y, \SAT_n(z,y)$. $A$ is $A^{\SAT_n,s}_n$ if $\SAT_n\notin \Circuit[s^3]$ and an arbitrary $poly(s)$-size subset of $\{0,1\}^n$ otherwise. $C(z)$ in the right disjunct of $\alpha^s_n$ stands for the output of a single-output $s$-size circuit $C$ with input $z$. The circuit $C$ is represented by free variables. The circuit $B$ in the left disjunct is a fixed $poly(s)$-size circuit, with $n$ inputs and $\le n$ outputs, obtained from a fixed $s^3$-size single-output circuit $B'$ with $n$ inputs such that $$\SAT_n\in\Circuit[s^3]\Leftrightarrow \forall x\in\{0,1\}^n, B'(x)=\SAT_n(x).$$ The circuit $B'$ exists by considering two cases: I. If $\SAT_n\notin\Circuit[s^3]$, we can take arbitrary $s^3$-size circuit $B'$; II. If $\SAT_n\in \Circuit[s^3]$, we can let $B'$ be an $s^3$-size circuit computing $\SAT_n$. $B$ is obtained from $B'$ in a standard way so that $B$ solves the search version of $\SAT_n$ if $\SAT_n\in \Circuit[s^3]$. Similarly as before, we get the following.

\begin{theorem}[Circuit complexity from nonuniform proof complexity]\label{t:acore}\ \newline
Let $k\ge 3$ be a constant. If there are tautologies without p-size \EF-derivations from substitutional instances of tautologies $\alpha^{n^k}_n$, then $\SAT_n\notin\Circuit[n^{k}]$ for infinitely many $n$.
\end{theorem}

\noindent Theorem \ref{t:acore} shows that $\EF+\alpha^{n^k}$, defined analogously as $\EF+w^k(f)$, is in certain sense optimal: If $\forall n, \SAT_n\in\Circuit[n^k]$, then $\EF+\alpha^{n^k}$ has $poly(n)$-size proofs of all tautologies. In fact, there is a p-size circuit which given a tautology $\phi$ of size $n$ outputs its proof in $\EF+\alpha^{n^k}$. Cook and Kraj\'i\v{c}ek \cite{CK} constructed an optimal proof system with 1 bit of nonuniform advice. Their system differs from $\EF+\alpha^{n^k}$ in that it is based on a diagonalization simulating all possible proof systems.

\subsubsection{Meta-mathematical implications}

We use observations from the previous section to shed light on some classical questions in complexity theory. Specifically, we show that the assumption of the existence of a p-time witnessing function in Theorem \ref{t:core} can be replaced by an efficient algorithm generating anticheckers for \SAT or by efficient reductions collapsing some of Impagliazzo's worlds (together with a hardness of {\sf E} for subexponential-size circuits).
\bigskip

\noindent {\bf I. Feasible anticheckers.} 

\begin{theorem}[`CC $\leftarrow$ PC' from feasible anticheckers - Informal, cf. Theorem \ref{t:fant}]\label{t:ifant}\ \newline
Let $k\ge 3$ be a constant and assume that there is a p-time function $f$ such that \SB proves:
\begin{itemize}
\item[] ``$\forall 1^n$, $f(1^n)$ outputs a $poly(n^{k})$-size circuit $B$ such that $$\forall x,y \in \{0,1\}^n, [\SAT_n(x,y)\rightarrow \SAT_n(x,B(x))]$$ or $\big(f(1^n)$ outputs sets $A^{\SAT_n,n^k}_n,A'\subseteq \{0,1\}^n$ of size $poly(n^k)$ such that $A'$ is a set of assignment $y_x$ of formulas $x\in A^{\SAT_n,n^k}$ such that $y_x$ satisfies $x$ if $x$ is satisfiable, and
$\forall\ n^{k}$-size circuit $C$, $$\exists x\in A^{\SAT_n,n^k}_n, \SAT_n(x,y_x)\ne C(x)\big)."$$
\end{itemize}
Then, proving that \EF is not p-bounded implies $\SAT_n\notin\Circuit[n^{k}]$ for infinitely many $n$.
\end{theorem}

\def\prfifant{
\proof
The statement assumed to have an intuitionistic  \SB-proof is $\forall\Sigma^b_2$. By the witnessing theorem for intuitionistic \SB \cite{}, there is a p-time function $f$, such that intuitionistic \SB proves:
``$\forall 1^n>n_0$, $f(1^n)$ outputs a $poly(n^{3k})$-size circuit $B$ and $A^{\SAT_n,n^k}_n,A'\subseteq \{0,1\}^n, D\subseteq A^{\SAT_n,n^k}_n\times A'$ of size $poly(n)$ such that  $\forall x,y\in \{0,1\}^n, \forall n^k$-size circuit $C$,
$$\big(\SAT_n(x,y)\rightarrow \SAT_n(x,B(x))\big)\vee \big(D'\wedge \exists z\in A^{\SAT_n,n^k}_n, \SAT_n(z,y_z)\ne C(z)\big)."$$ Here, $D'$ stands for the properties of $D$ guaranteed by the assumption. Since this is a $\Pi^b_1$-statement, there are p-size \EF-proofs of its propositional translation.

If we now assume that $\forall n, \SAT_n\in \Circuit [n^k]$, there is a circuit $C$ falsifying the right disjunct of the displayed formula (or $y$ falsifying $D'$). Consequently, \EF proves efficiently that the circuit generated by $f(1^n)$ solves $\SAT_n$, which implies that \EF is p-bounded.
\qed

\medskip}

Theorem \ref{t:ifant} follows easily from the proof of Theorem \ref{t:core}, see \S\ref{s:fant}. 

Similarly as in Corollary \ref{c:clbfromwit}, if for each $k\ge 3$ there is a p-time function $f$ such that \SB proves the statement about the existence of anticheckers from the assumption of Theorem \ref{t:ifant}, then $\NP\not\subseteq {\sf SIZE}[n^{\ell}]$, for all $\ell\in\mathbb{N}$.

\def\fantinpv{
If we used \PV instead of intuitionistic \SB in Theorem \ref{t:fant}, by the KPT theorem \ref{}, the \PV-provability of the existence of anticheckers would imply the existence of p-time functions $f_1,\dots,f_{c}$, for a constant $c$, with a \PV-proof of:
\begin{itemize}
\item[] ``$\forall 1^n>n_0$, $\forall x^1,\dots,x^c,y^1,\dots,y^c \in \{0,1\}^n, \forall n^k$-size circuits $C^1,\dots, C^c$, $f_1(1^n)$ outputs a $poly(n^{3k})$-size circuit $B$ and $A^{\SAT_n,n^k}_n,A'\subseteq \{0,1\}^n, D\subseteq A^{\SAT_n,n^k}_n\times A'$ of size $poly(n)$ such that $P_{f_1}(x^1,y^1,C^1)$ holds, where $$
P_{f_1}(x,y,C):=\big(\SAT_n(x,y)\rightarrow \SAT_n(x,B(x))\big)\vee \big(D'\wedge \exists z\in A^{\SAT_n,n^k}_n, \SAT_n(z,y_z)\ne C(z)\big),$$ or $f_2(1^n, x^1,y^1,C^1)$ outputs a $poly(n^{3k})$-size circuit $B$ and $A^{\SAT_n,n^k}_n,A',D$ of size $poly(n)$ such that $P_{f_2}(x^2,y^2,C^2)$ holds, or \\
$\dots$\\
or $f_c(1^n,x^1,\dots,x^{c-1},y^1,\dots,y^{c-1},C^1,\dots,C^{c-1})$ outputs a $poly(n^{3k})$-size circuit $B$ and $A^{\SAT_n,n^k}_n,A', D$ of size $poly(n)$ such that $P_{f_c}(x^c,y^c,C^c)$."
\end{itemize}
The resulting $\Pi^b_1$-statement could be translated to propositional tautologies with p-size \EF-proofs. 
However, given $\forall n, \SAT_n\in \Circuit [n^k]$, we could not directly obtain p-size \EF-proofs of tautologies stating that one of the functions $f_1,\dots,f_c$ generates a circuit solving $\SAT_n$. This is because $B$ and $A^{\SAT_n,n^k}_n$ generated by $f_2$ depend on $y^1$. For the same reason, it seems possible for \EF to prove efficiently $\SAT\in \Ppoly$ (using the formalization based on the `KPT witnessing') without proving efficiently all tautologies. (formulating anticheckers as 'ex A st forall C C solves sat or C fails on A'  deals with 1st disjunct)}
\bigskip

\noindent {\bf II. One-way functions from $\NP\not\subseteq\Ppoly$.}
\medskip

\noindent Denote by $\ttable(f_n,s)$ a propositional formula expressing that Boolean function $f_n:\{0,1\}^n\mapsto \{0,1\}$ represented by its truth-table is not computable by a Boolean circuit of size $s$ represented by free variables, see \S $\ref{s:formalizations}$. So $\ttable(f_n,s)$ is a tautology if and only if $f_n$ is hard for circuits of size $s$. The size of the formula $\ttable(f_n,s)$ is $poly(2^n,s)$. Similarly, let $\ttable(f_n,s,t)$ be a formula expressing that circuits of size $s$ fail to compute $f$ on $\ge t$-fraction of inputs.

Given a function $h\in {\sf E}$ such that for some $n_0$, for each $n\ge n_0$, each $s(2^n)$-size circuit with $n$ inputs fails to compute $h_n$ on $\ge t(2^n)$-fraction of inputs, where $s,t$ are p-time functions in $2^n$, we define a proof system $\EF+\ttable(h,s,t)$ as an extension of \EF which is allowed to derive in its proofs substitutional instances of $\ttable(h_n,s(2^n),t(2^n))$, for $n\ge n_0$.

\begin{theorem}[`CC $\leftarrow$ PC' from `OWF $\leftarrow$ $\NP\not\subseteq \Ppoly$' \& hard {\sf E}, cf. Theorem \ref{t:owpnp}]\label{t:iowpnp}\ \newline
Assume that for each sufficiently big $n$, each $2^{n/4}$-size circuit fails to compute $h'\in {\sf E}$ on $\ge 1/2-1/2^{n/4}$ of all inputs. Further, assume that there is a p-time function $h:\{0,1\}^n\mapsto \{0,1\}^{u(n)}$ such that for each constants $c,d$, there is a p-time function $f$ and a constant $0<\epsilon<1$ such that \SB proves:
\begin{itemize}
\item[] ``$\forall n, \forall cn^{c}$-size circuit $C$ with $u(m)$ inputs and $m$ outputs such that $n\le dm^d$,\\
$f(C)$ outputs a $poly(n)$-size circuit $B$ such that $$\forall x,y \in \{0,1\}^{n^{\epsilon}}, [\SAT_{n^{\epsilon}}(x,y)\rightarrow \SAT_{n^{\epsilon}}(x,B(x))]$$ or $$\Pr_{x\in \{0,1\}^{m}}[h(C(h(x)))=h(x)]<1/2."$$
\end{itemize}
Then, proving that $\EF+\ttable(h',2^{n/4},1/2-1/2^{n/4})$ is not p-bounded implies $\SAT\notin\Ppoly$.
\end{theorem}

We remark that if there is a p-time function $f'$ and constant $n_0$ such that \SB proves that for each $n\ge n_0, f'(1^{2^n})$ outputs the truth-table of a function $h'$ with $n$ inputs which is hard on average for $2^{n/4}$-size circuits, then \EF and $\EF+\ttable(h',2^{n/4},1/2-1/2^{n/4})$ are p-equivalent. (Note that $h'\in \mathsf{E}$.)

The proof of Theorem \ref{t:iowpnp} is based on a formalization of the already-mentioned fact that given a one-way function $h$ and a function in {\sf E} hard for subexponential-size circuits, we can construct a p-time function witnessing errors of p-size circuits attempting to solve \SAT. (The witnessing function outputs formulas encoding the statement $h(x)=b$, with free variables $x$ and suitable constants $b$, cf. Lemma \ref{l:owpwit}.) We formalize the conditional witnessing in a theory $\SB+dWPHP(\pv)$, where $dWPHP(\pv)$ stands for a dual weak pigeonhole principle, see \S\ref{sba}. Combining this with the assumption that \SB proves that a one-way function can be obtained from the hardness of \SAT, we obtain the `ideal' $(\SB+dWPHP(\pv))$-provable witnessing similar to the tautology $w^k_n(f)$. Having the ideal witnessing statement we proceed as in the proof of Theorem \ref{t:core} with the difference that the axiom $dWPHP(\pv)$ and the assumed hardness of ${\sf E}$ lead to the system $\EF+\ttable(h',2^{n/4},1/2-1/2^{n/4})$ instead of \EF, see \S\ref{s:owpnp}. 

\bigskip
\noindent {\bf III. Learning from the non-existence of OWFs.}


\begin{theorem}[`CC $\leftarrow$ PC' from `Learning $\leftarrow$ $\not\exists$ OWF' \& hardness of E, cf. Theorem \ref{t:dich}]\label{t:idich}
Let $k,t\ge 1$ be constants. Assume that for each sufficiently big $n$, each $2^{n/4}$-size circuit fails to compute $h'\in {\sf E}$ on $\ge 1/2-1/2^{n/4}$ of all inputs. Further, assume that there is a p-time functions $h:\{0,1\}^n\mapsto \{0,1\}^{u(n)}$ such that for each constants $c,d$, there is a p-time function $f$ and a constant $0<\epsilon<1$ such that \SB proves:
\begin{itemize}
\item[] ``$\forall n$, $\forall cn^c$-size circuits $C$ with $u(m)$ inputs and $m$ outputs such that $n\le dm^d$,\\ $f(C)$ outputs a $poly(n)$-size circuit $B$ such that $B$ learns $n^{\epsilon t}$-size circuits with $n^{\epsilon}$ inputs, over the uniform distribution, up to error $1/2-1/n^{\epsilon}$, with membership queries and confidence $1/n^{\epsilon}$, \\ or $$\Pr_{x\in \{0,1\}^m}[h(C(h(x)))=h(x)]<1/2."$$
\end{itemize}
Then, there are constants $b$ and $a$ (depending on $k,t,h,h',c,d,f,\epsilon$) such that for each $n$, the existence of a function $g_n:\{0,1\}^n\mapsto \{0,1\}$ such that no $bn^b$-size circuit computes $g_n$ on $\ge (1/2+1/n)$ fraction of inputs and such that $\EF+\ttable(h',2^{n/4},1/2-1/2^{n/4})$ does not have a $2^{an}$-size proof of $\ttable(g_n,n^{t})$ implies that 
$\SAT_n\notin\Circuit[n^k]$.
\end{theorem}

Theorem \ref{t:idich} is proved similarly as Theorem \ref{t:iowpnp} with the difference that the provability of efficient learning allows us to prove efficiently only circuit lower bounds instead of all tautologies, see \S\ref{s:dich}.

\cite[Lemma 4]{PSlearaut} shows that, assuming {\sf E} is \SB-provably hard as in Theorem \ref{t:idich}, learning algorithms for small circuits can be \SB-provably constructed from circuits automating \EF on $\ttable$-formulas.\footnote{\cite[Lemma 4]{PSlearaut} assumes also the existence of a prime. The assumption can be removed after moving to the propositional setting.} Theorem \ref{t:idich} thus implies that (assuming ${\sf E}$ is \SB-provably hard) \SB-deriving automatability of \EF from the non-existence of one-way functions would reduce circuit complexity to \EF lower bounds.

\bigskip

\noindent{\it Generalization to stronger proof systems.} \SB in Theorems \ref{t:core} \& \ref{t:ifant}-\ref{t:idich} can be replaced by essentially an arbitrary first-order theory $T$ containing \SB and satisfying some basic properties, if we simultaneously replace \EF in conclusions of Theorems \ref{t:core} \& \ref{t:ifant}-\ref{t:idich} by a suitable propositional proof system $P_T$ such that propositional translations of $\Pi^b_1$ theorems of $T$ have p-size proofs in $P_T$. 
\medskip

\noindent {\it Weakening the assumptions.} The core component of Theorems \ref{t:core} \& \ref{t:ifant}-\ref{t:idich} is the existence of a suitable reduction. For example, in case of Theorem \ref{t:idich} we need a p-time reduction constructing learning algorithms from circuits breaking one-way functions. If such a reduction exists, even without assuming its provability in \SB, we can build a propositional proof system $P$ by adding tautologies encoding the correctness of the reduction to \EF. 
Then, showing that the resulting proof system $P+\ttable(h',2^{n/4},1/2-1/2^{n/4})$ is not p-bounded on $\ttable$-tautologies would separate \Ptime and \NP. This shows that the \SB-provability in the assumptions of our theorems can be weakened just to the validity of the respective statements, if we use stronger systems than \EF in their conclusions. It also shows that the most of technicalities in the present paper stem from making the presented approach work for \EF.

Moreover, we remark that if our final goal is to prove that $\Ptime\ne\NP$, then the first assumption of Theorems \ref{t:iowpnp}-\ref{t:idich} postulating the existence of a hard Boolean function in {\sf E} is given to us `for free', as otherwise, if all functions in {\sf E} can be approximated by subexponential-size circuits, it is not hard to show that $\Ptime\ne\NP$.

\subsubsection{Plausibility of the assumptions}\label{s:plausibility}

\noindent {\bf Impagliazzo's worlds.} In a famous survey of Impagliazzo \cite{Iw}, he described 5 possible worlds of average-case complexity: Algorithmica, Heuristica, Pessiland, Minicrypt and Cryptomania. Recently, there have been various approaches proposed to rule out Heuristica and Pessiland (see, for example, \cite{Hb,S19,LP}) by studying the complexity of problems about compression, such as the Minimum Circuit Size Problem (MCSP) and the problem of computing time-bounded Kolmogorov complexity. Our results have implications for the feasibility of such efforts - provable collapses of Impagliazzo's worlds would imply a new and surprising link between proof complexity and circuit complexity. For example, the reduction assumed in Theorem \ref{t:idich} asks for a construction of learning algorithms from circuits breaking one-way functions. Morally, the existence of such a reduction would rule out Pessiland out of Impagliazzo's worlds. Similarly, the reduction assumed in Theorem \ref{t:iowpnp} would rule out Pessiland and Heuristica.
\bigskip

\noindent {\bf Feasible MinMax theorem.} The proofs of the existence of anticheckers we are aware of use the efficient MinMax theorem \cite{Alt,LY,Nm} or similar methods. If we had a proof of MinMax which would use counting with only polynomial precision (formally, \APC-counting) and if we could replace p-time $f$ in Theorem \ref{t:ifant} by the existential quantifiers (see \S\ref{s:fant} for a discussion of the issue), we could prove the existence of anticheckers in \APC and obtain the desired reduction of circuit complexity to proof complexity. Here, \APC is Je\v{r}\'abek's theory of approximate counting, see \S\ref{sba}.
\bigskip

\def\previntro{
\noindent {\bf Learning algorithms.} In the PAC model of learning introduced by Valiant \cite{LT}, a circuit class $\mathcal{C}$ is learnable by a randomized algorithm $L$ over the uniform distribution, up to error $\epsilon$, with confidence $\delta$ and membership queries, if for every Boolean function $f$ computable by a circuit from $\mathcal{C}$, when given oracle access to $f$, $L$ outputs with probability $\ge \delta$ over the uniform distribution a circuit computing $f$ on $\ge (1-\epsilon)$ inputs. An important task of learning theory is to find out if standard circuit classes such as \Ppoly  are learnable by efficient circuits. A way to approach the question is to connect the existence of efficient learning algorithms to other standard conjectures in complexity theory. For example, we can try to prove that efficient learning of \Ppoly is equivalent to $\Ptime =\NP$ or to the non-existence of strong pseudorandom generators. In both cases one implication is known: $\Ptime=\NP$ implies efficient learning of \Ppoly (with small error and high confidence)
 which in turn breaks pseudorandom generators. However, while some progress on the opposite implications has been made, they remain open, cf. \cite{ABX,S19}. 
}

\noindent {\bf \SB-provability of a circuit lower bound.} If we want to replace $\EF+\ttable(h',2^{n/4},1/2-1/2^{n/4})$ in Theorems \ref{t:iowpnp} \& \ref{t:idich} by \EF, it suffices to assume the \SB-provability of a subexponential circuit lower bound for {\sf E}. This assumption has an interesting status. Razborov's conjecture about hardness of Nisan-Wigderson generators implies a conditional hardness of formulas $\ttable(h,n^{O(1)})$ for Frege (for every $h$), cf. \cite{Rkdnf}, and it is possible to consider extensions of the conjecture to all standard proof systems, even set theory $\mathsf{ZFC}$. A conditional hardness of $\ttable$-formulas (for \EF) follows also from a conjecture of Kraj\'i\v{c}ek \cite[Conjecture 7.9]{Kwphp}. If the $\ttable$-formulas expressing subexponential lower bounds for ${\sf E}$ are hard for \EF, then \SB cannot prove the lower bounds either. On the other hand, it is not known how to prove hardness of $\ttable(h,2^{n/4})$, 
 for all $h$, for Frege, under any standard complexity-theoretic hardness assumption. Moreover, all major circuit lower bounds for weak circuit classes and explicit Boolean functions are known to be provable in \SB, cf. \cite{Rup,MP}.\footnote{This has not been verified for lower bounds obtained via the algorithmic method of Williams \cite{Wacc}.} It is thus perfectly possible that subexponential average-case circuit lower bounds for $\mathsf{E}$ are provable in a theory such as \SB.\footnote{We emphasize that Theorems \ref{t:iowpnp} \& \ref{t:idich} do not require that \SB proves a circuit lower bound. Further, we can replace the system $\EF+\ttable(h',2^{n/4},1/2-1/2^{n/4})$ in Theorems \ref{t:iowpnp} \& \ref{t:idich} by, say, a propositional proof system corresponding to {\sf ZFC}, if we assume that a subexponential circuit lower bound for {\sf E} is provable in (the set theory) {\sf ZFC}.}

\subsubsection{Revising the status of the Cook-Reckhow program}\label{ss:CRstatus}

Showing that statements like $\Ptime\ne\NP$ follow from proof complexity lower bounds for concrete proof systems is considered so challenging that there have not been practically any conscious attempts to approach it. Our results show that the significant efforts that have been made in order to address some of the central problems in cryptography and learning theory are, in fact, aiming to establish precisely that. This can be interpreted as an evidence for the hardness of resolving the relevant problems in cryptography and learning theory, but also as showing that proving that the Cook-Reckhow program could be in principle realized successfully might not be completely out of reach. In any case, the presented results demonstrate a new fundamental connection between proof complexity, cryptography and learning theory.

\subsubsection{Self-provability from random self-reducibility}

In \S\ref{s:rselfred} we show that a self-provability of circuit upper bounds can be obtained for some random-selfreducible problems. We exhibit this in the case of the discrete logarithm: Assuming the existence of a hard Boolean function in {\sf E} and the existence of efficient circuits solving the discrete logarithm problem, we show that an explicit proof system admits p-size proofs of the fact that the discrete logarithm problem is solvable efficiently.

\subsection{Open problems} 1. Prove that any superpolynomial \EF lower bound separates \Ptime and \NP by constructing a witnessing function such that, for all big enough $n$, tautologies $w^k_n(f)$ have p-size \EF-proofs. It would be interesting to construct such witnessing functions $f$ even assuming $\Ptime\not\subseteq {\sf SIZE}[n^{\ell}]$, for all $\ell\in\mathbb{N}$. By Corollary \ref{c:clbfromwit}, an unconditional construction of such witnessing functions would imply $\NP\not\subseteq {\sf SIZE}[n^{\ell}]$, for all $\ell\in\mathbb{N}$.
\smallskip

\noindent 2. Prove that the p-time functions defining reductions in Theorems \ref{t:ifant}-\ref{t:idich} can be replaced by the existential quantifiers, see \S\ref{s:fant} for a discussion of the problem.

\section{Preliminaries}

\subsection{Learning algorithms}\label{ss:npl}

$[n]$ denotes $\{1,\dots,n\}$. $1^n$ stands for a string of $n$ 1s. $\Circuit[s]$ denotes fan-in two Boolean circuits of size at most $s$. The size of a circuit is the number of its gates. A function $f:\{0,1\}^n\mapsto \{0,1\}$ is $\gamma$-approximated by a circuit $C$, if $\Pr_x[C(x)=f(x)]\ge\gamma$.

\def\natprfdef{
\begin{definition}[Natural property \cite{RR}]\label{d:natpr} Let $m=2^n$ and $s,d:\mathbb{N} \mapsto \mathbb{N}$. A sequence of circuits $\{C_{m}\}^{\infty}_{n=1}$ is a $\Circuit[s(m)]$-natural property useful against $\Circuit[d(n)]$ if
\begin{itemize}
\item[1.] {\em Constructivity.} $C_{m}$ has $m$ inputs and size $s(m)$,
\item[2.] {\em Largeness.} $\Pr_x[C_{m}(x)=1]\ge 1/m^{O(1)}$,
\item[3.] {\em Usefulness.} For each sufficiently big $m$, $C_{m}(x)=1$ implies that $x$ is a truth-table of a function on $n$ inputs which is not computable by circuits of size $d(n)$.
\end{itemize}
\end{definition}
}


\begin{definition}[PAC learning]\label{d:lear}
A circuit class $\mathcal{C}$ is learnable over the uniform distribution by a circuit class $\mathcal{D}$ up to error $\epsilon$ with confidence $\delta$, if there are randomized oracle circuits $L^f$ from $\mathcal{D}$ such that for every Boolean function $f:\{0,1\}^n\mapsto\{0,1\}$ computable by a circuit from $\mathcal{C}$, when given oracle access to $f$, input $1^n$ and the internal randomness $w \in \{0,1\}^*$, $L^f$ outputs the description of a circuit satisfying 
  \begin{equation*}
    \Pr_w [ L^f(1^n, w) \text{ } (1-\epsilon) \text{-approximates } f ] \geq \delta.
  \end{equation*}

\noindent $L^f$ uses non-adaptive membership queries if the set of queries which $L^f$ makes to the oracle does not depend on the answers to previous queries. $L^f$ uses random examples if the set of queries which $L^f$ makes to the oracle is chosen uniformly at random. 
\end{definition}

In this paper, PAC learning always refers to learning over the uniform distribution. While, a priori, learning over the uniform distribution might not reflect real-world scenarios very well (and on the opposite end, learning over all distributions is perhaps overly restrictive), as far as we can tell it is possible that PAC learning of p-size circuits over the uniform distribution implies PAC learning of p-size circuits over all p-samplable distributions. Binnendyk, Carmosino, Kolokolova, Ramyaa and Sabin \cite{BCKRS} proved the implication, if the learning algorithm in the conclusion is allowed to depend on the p-samplable distribution.

\def\boosting{
\noindent {\bf Boosting confidence and reducing error.} The confidence of the learner can be efficiently boosted in a standard way. Suppose an $s$-size circuit $L^f$ learns $f$ up to error $\epsilon$ with confidence $\delta$. We can then run $L^f$ $k$ times, test the output of $L^f$ from every run with $m$ new random queries and output the most accurate one. By Hoeffding's inequality, $m$ random queries fail to estimate the error $\epsilon$ of an output of $L^f$ up to $\gamma$ with probability at most $2/e^{2\gamma^2m}$. 
Therefore the resulting circuit of size $poly(s,m,k)$ learns $f$ up to error $\epsilon+\gamma$ with confidence at least $1-2k/e^{2\gamma^2m}-(1-\delta)^k\ge 1-2k/e^{2\gamma^2m}-e^{-k\delta}$. If we are trying to learn small circuits  we can get even confidence 1 by fixing internal randomness of learner nonuniformly without losing much on the running time or the error of the output. It is also possible to reduce the error up to which $L^f$ learns $f$ without a significant blowup in the running time and confidence. If we want to learn $f$ with a better error, we first learn an amplified version of $f$, $Amp(f)$. Employing direct product theorems and Goldreich-Levin reconstruction algorithm, Carmosino et. al. \cite[Lemma 3.5]{CIKK} showed that for each $0<\epsilon,\gamma<1$ it is possible to map a Boolean function $f$ with $n$ inputs to a Boolean function $Amp(f)$ with $poly(n,1/\epsilon,\log(1/\gamma))$ inputs so that $Amp(f)\in \Ppoly^f$ and there is a probabilistic $poly(|C|,n,1/\epsilon,1/\gamma)$-time machine which given a circuit $C$ $(1/2+\gamma)$-approximating $Amp(f)$ and an oracle access to $f$ outputs with high probability a circuit $(1-\epsilon)$-approximating $f$. 
We can thus often ignore the optimisation of the confidence and error parameter. Note, however, that the error reduction of Carmosino et al. requires membership queries.
}

\def\cikkref{
\bigskip

\noindent{\bf Natural proofs vs learning algorithms.}
Natural proofs are actually equivalent to efficient learning algorithms with suitable parameters. In this paper we need just one implication.

\begin{theorem}[Carmosino-Impagliazzo-Kabanets-Kolokolova \cite{CIKK}]\label{t:cikkcore}
Let $R$ be a $\Ppoly$-natural property useful against $\Circuit[n^k]$ for $k \geq 1$. Then, for each $\gamma\in (0,1)$, $\Circuit[n^{k\gamma/a}]$ is learnable by $\Circuit[2^{O(n^{\gamma})}]$ over the uniform distribution with non-adaptive membership queries, confidence 1, up to error $1/n^{k\gamma/a}$, where 
$a$ is an absolute constant.
\end{theorem}
}

\subsection{Bounded arithmetic and propositional logic}\label{sba}

Theories of bounded arithmetic capture various levels of feasible reasoning and present a uniform counterpart to propositional proof systems.

The first theory formalizing p-time reasoning was introduced by Cook \cite{Cpv} as an equational theory \pv. We work with its first-order conservative extension \PV from \cite{KPT}. The language of \PV, denoted \pv as well, consists of symbols for all p-time algorithms given by Cobham's characterization of p-time functions, cf. \cite{Cptime}. A \pv-formula is a first-order formula in the language \pv. $\Sigma^b_0$ (=$\Pi^b_0$) denotes \pv-formulas with only sharply bounded quantifiers $\exists x, x\leq |t|$, $\forall x, x\leq |t|$, where 
$|t|$ is ``the length of the binary representation of $t$". Inductively, $\Sigma^b_{i+1}$ resp. $\Pi^b_{i+1}$ is the closure of $\Pi^b_i$ resp. $\Sigma^b_i$ under positive Boolean combinations, sharply bounded quantifiers, and bounded quantifiers $\exists x, x\le t$ resp. $\forall x, x\le t$. 
Predicates definable by $\Sigma^b_{i}$ resp. $\Pi^b_{i}$ formulas are in the $\Sigma^p_{i}$ resp. $\Pi^p_{i}$ level of the polynomial hierarchy, and vice versa. 
\PV is known to prove $\Sigma^b_0(\pv)$-induction: $$A(0)\wedge \forall x\ (A(x)\rightarrow A(x+1))\rightarrow \forall x A(x),$$ for $\Sigma^b_0$-formulas $A$, cf. Kraj\'{i}\v{c}ek \cite{Kba}. 

Buss \cite{Bba} introduced the theory \SB extending \PV with the $\Sigma^b_1$-length induction: $$A(0)\wedge \forall x<|a|, (A(x)\rightarrow A(x+1))\rightarrow A(|a|),$$ for $A\in\Sigma^b_1$. \SB proves the sharply bounded collection scheme $BB(\Sigma^b_1)$: $$\forall i<|a|\ \exists x<a, A(i,x)\rightarrow \exists w\ \forall i<|a|, A(i,[w]_i),$$ for $A\in\Sigma^b_1$ ($[w]_i$ is the $i$th element of the sequence coded by $w$), which is unprovable in \PV under a cryptographic assumption, cf. \cite{CTcol}. On the other hand, \SB is $\forall\Sigma^b_1$-conservative over \PV. This is a consequence of Buss's witnessing theorem stating that $\SB\vdash \exists y, A(x,y)$ for $A\in\Sigma^b_1$ implies $\PV\vdash A(x,f(x))$ for some \pv-function $f$. 

Following a work by Kraj\'{i}\v{c}ek \cite{Kwphp}, Je\v{r}\'{a}bek \cite{Jwphp,Jphd,Japx} systematically developed a theory \APC capturing probabilistic p-time reasoning by means of approximate counting.\footnote{Kraj\'{i}\v{c}ek \cite{Kwphp} introduced a theory $BT$ defined as $\SB+dWPHP(\pv)$ and proposed it as a theory for probabilistic p-time reasoning.} The theory \APC is defined as $\PV+dWPHP(\pv)$ where $dWPHP(\pv)$ stands for the dual (surjective) pigeonhole principle for \pv-functions, i.e. for the set of all formulas $$x>0\rightarrow\exists v<x(|y|+1)\forall u<x|y|,\ f(u)\neq v,$$ where $f$ is a \pv-function which might involve other parameters not explicitly shown. We devote \S\ref{approximatec} to a more detailed description of the machinery of approximate counting in \APC.
\bigskip

Any $\Pi^b_1$-formula $\phi$ provable in \SB can be expressed as a sequence of tautologies $||\phi||_n$ with proofs in the Extended Frege system \EF which are constructible in p-time (given a string of the length $n$), cf. \cite{Cpv}. We refer to Kraj\'i\v{c}ek \cite{Kpc} for basic notions in proof complexity such as \EF.
 \def\WFstuff{Similarly, $\Pi^b_1$-formulas provable in \APC translate to tautologies with p-time constructible proofs in \WF, an extension of \EF introduced by Je\v{r}\'{a}bek \cite{Jwphp}. We describe the translation and system \WF in more detail below. 
}
As it is often easier to present a proof in a theory of bounded arithmetic than in the corresponding propositional system, bounded arithmetic functions, so to speak, as a uniform language for propositional logic.
\def\WFstuff2{\smallskip

We refer to Kraj\'i\v{c}ek \cite{Kpc} for basic notions in proof complexity.

The translation of a $\Pi^b_1$ formula $\phi$ into a sequence of propositional formulas $||\phi||_{\overline{n}}$ works as follows (full details can be found in \cite[Section 12.3]{Kpc}). For each $\pv$-function $f(x_1,\dots,x_k)$ and numbers $n_1,\dots,n_k$ we have a p-size circuit $C_f$ computing the restriction $f:2^{n_1}\times \dots\times 2^{n_k}\mapsto 2^{b(n_1,\dots,n_k)}$, where $b$ is a suitable `bounding' polynomial for $f$. The formula $||f||_{\overline{n}}(p,q,r)$ expresses that $C_f$ outputs $r$ on input $p$, with $q$ being the auxiliary variables corresponding to the nodes of $C_f$.  The formula $||\phi(x)||_{\overline{n}}(p,q)$ is defined as $||\phi'(x)||_{\overline{n}}(p,q)$, where $\phi'$ is the negation normal form of $\phi$, i.e. negations in $\phi'$ are only in front of atomic formulas. The formula $||\phi'(x)||_{\overline{n}}(p,q)$ is defined inductively in a straightforward way so that $||\dots||$ commutes with $\vee,\wedge$. The atoms $p$ correspond to variables $x$, atoms $q$ correspond to the universally quantified variables of $\phi$ and to the outputs and auxiliary variables of circuits $C_f$ for functions $f$ appearing in $\phi$. Sharply bounded quantifiers are replaced by polynomially big conjuctions resp. disjunctions. For the atomic formulas we have\footnote{Formally, we should allow terms, not just function symbols $f,g$, in atomic formulas.}, \begin{align*}||f(x)=g(x)||_{\overline{n}}:= & ||f(x)||_{\overline{n}}(p,q,r)\wedge ||g(x)||_{\overline{n}}(p,q',r')\rightarrow \bigwedge_i r_i=r'_i,\\
||\neg f(x)=g(x)||_{\overline{n}}:= & ||f(x)||_{\overline{n}}(p,q,r)\wedge ||g(x)||_{\overline{n}}(p,q',r')\rightarrow \neg\bigwedge_i r_i=r'_i,\\
||f(x)\le g(x)||_{\overline{n}}:= & ||f(x)||_{\overline{n}}(p,q,r)\wedge ||g(x)||_{\overline{n}}(p,q',r')\rightarrow \bigwedge_i (r_i\wedge \bigwedge_{j>i}(r_j=r_j')\rightarrow r'_i),\\
||\neg f(x)\le g(x)||_{\overline{n}}:= & ||f(x)||_{\overline{n}}(p,q,r)\wedge ||g(x)||_{\overline{n}}(p,q',r')\rightarrow \neg \bigwedge_i(r_i\wedge\bigwedge_{j>i}(r_j=r_j')\rightarrow r'_i).\end{align*}
}

\subsection{Approximate counting}\label{approximatec}
In order to prove our results 
we will need to use Je\v{r}\'{a}bek's theory of approximate counting. This section recalls the properties of \APC we will need.
\medskip



By a definable set we mean a collection of numbers satisfying some formula, possibly with parameters. When a number $a$ is used in a context which asks for a set it is assumed to represent the integer interval $[0,a)$, e.g. $X\subseteq a$ means that all elements of set $X$ are less than $a$. If $X\subseteq a$, $Y\subseteq b$, then $X\times Y:=\{bx+y\mid x\in X, y\in Y\}\subseteq ab$ and $X\dot{\cup} Y:=X\cup\{y+a\mid y\in Y\}\subseteq a+b$. Rational numbers are assumed to be represented by pairs of integers in the natural way. We use the notation $x\in Log\leftrightarrow \exists y,\ x=|y|$ and $x\in LogLog\leftrightarrow \exists y,\ x=||y||$.
\medskip

Let $C: 2^n\rightarrow 2^m$ be a circuit and $X\subseteq 2^n, Y\subseteq 2^m$ definable sets. 
We write $C:X\twoheadrightarrow Y$ if $\forall y\in Y\exists x\in X,\ C(x)=y$.
Je\v{r}\'{a}bek \cite{Japx} gives the following definitions in \APC, but they can be formulated in \PV as well.



\begin{definition} Let $X,Y\subseteq 2^n$ be definable sets, and $\epsilon\leq 1$. The size of $X$ is approximately less than the size of $Y$ with error $\epsilon$, written as $X\preceq_{\epsilon} Y$, if there exists a circuit $C$, and $v\neq 0$ such that $$C: v\times (Y\dot{\cup}\epsilon 2^n)\twoheadrightarrow v\times X.$$ $X\approx_{\epsilon}Y$ stands for $X\preceq_{\epsilon} Y$ and $Y\preceq_{\epsilon} X$.
\end{definition}

Since a number $s$ is identified with the interval $[0,s)$, $X\preceq_{\epsilon} s$ means that the size of $X$ is at most $s$ with error $\epsilon$.


The definition of $X\preceq_{\epsilon} Y$ is an unbounded $\exists \Pi^b_2$-formula even if $X,Y$ are defined by circuits so it cannot be used freely in bounded induction. Je\v{r}\'{a}bek \cite{Japx} solved this problem by working in ${\sf HARD}^A$, a conservative extension of \APC, defined as a relativized theory $\PV(\alpha)+dWPHP(\pv(\alpha))$ extended with axioms postulating that $\alpha(x)$ is a truth-table of a function on $||x||$ variables hard on average for circuits of size $2^{||x||/4}$, see \S\ref{ss:lform}. In ${\sf HARD}^A$ there is a $\PV(\alpha)$ function $Size$ approximating the size of any set $X\subseteq 2^n$ defined by a circuit $C$ so that $X\approx_{\epsilon} Size(C,2^n,2^{\epsilon^{-1}})$ for $\epsilon^{-1}\in Log$, cf. \cite[Lemma 2.14]{Japx}. If $X\cap t\subseteq 2^{|t|}$ is defined by a circuit $C$ and $\epsilon^{-1}\in Log$, we can define  $$\Pr_{x<t}[x\in X]_{\epsilon}:=\frac{1}{t} Size(C,2^{|t|},2^{\epsilon^{-1}}).$$

\def\relatth{
\begin{definition}[in \PV] Let $f:2^k\rightarrow 2$ be a truth-table of a Boolean function with $k$ inputs ($f$ is encoded as a string of $2^k$ bits, hence $k\in LogLog$). We say that $f$ is (worst-case) $\epsilon$-hard, written as $Hard_{\epsilon}(f)$ if no circuit $C$ of size $2^{\epsilon k}$ computes $f$. The function $f$ is average-case $\epsilon$-hard, written as $Hard^A_{\epsilon}(f)$, if for no circuit $C$ of size $\leq 2^{\epsilon k}$: $$|\{u<2^k|C(u)=f(u)\}|\geq (1/2+2^{-\epsilon k})2^k.$$
\end{definition}

\begin{proposition}[Je\v{r}\'{a}bek \cite{Jdual}] For every constant $\epsilon <1/3$ there exists a constant $c$ such that \APC proves: for every $k\in LogLog$ such that $k\geq c$, there exist average-case $\epsilon$-hard functions $f:2^k\rightarrow 2$.  
\end{proposition}

\PV can be relativized to $\PV(\alpha)$. The new function symbol $\alpha$ is then allowed in the inductive clauses for introduction of new function symbols. 
This means that the language of $\PV(\alpha)$, denoted also $\pv(\alpha)$, contains symbols for all p-time oracle algorithms.

\begin{definition}[Je\v{r}\'{a}bek \cite{Jwphp}] The theory ${\sf HARD}^A$ is an extension of the theory $\PV(\alpha)+dWPHP(\pv(\alpha))$ by the axioms
\smallskip

1. $\alpha(x)$ is a truth-table of a Boolean function in $||x||$ variables,
 
2. $x\geq c\rightarrow Hard^A_{1/4}(\alpha(x))$,

3. $||x||=||y||\rightarrow \alpha(x)=\alpha(y)$,
\smallskip

\noindent where $c$ is the constant from the previous proposition. 
\end{definition}

\begin{theorem}[Je\v{r}\'{a}bek \cite{Jdual,Japx}]\label{sec} ${\sf HARD}^A$ is a conservative extension of \APC. Moreover, there is a $\pv(\alpha)$-function $Size$ such that ${\sf HARD}^A$ proves: if $X\subseteq 2^n$ is definable by a circuit $C$, then $$X\approx_{\epsilon} Size(C,2^n,e)$$ where $\epsilon=|e|^{-1}$
\end{theorem} 

We will abuse the notation and write $Size(X,\epsilon)$ instead of $Size(C,2^n,e)$. 

\begin{definition}[in \APC] If $X\subseteq 2^{|t|}$ is defined by a circuit $C$ and $\epsilon^{-1}\in Log$, we put $$\Pr_{x<t}[x\in X]_{\epsilon}:=\frac{1}{t} Size(C,2^{|t|},\epsilon).$$ 
\end{definition}
}

The presented definitions of approximate counting are well-behaved:

\begin{proposition}[Je\v{r}\'{a}bek \cite{Japx}](in \PV)\label{lem} Let $X,X',Y,Y',Z\subseteq 2^n$ and $W,W'\subseteq 2^m$ be definable sets, and $\epsilon, \delta<1$. Then

$i)\ \ X\subseteq Y\Rightarrow X\preceq_{0} Y$,

$ii)\ \ X\preceq_{\epsilon} Y \wedge Y\preceq_{\delta} Z\Rightarrow X\preceq_{\epsilon+\delta} Z$,

$iii)\ \ X\preceq_\epsilon X'\wedge W\preceq_{\delta}W'\Rightarrow X\times W\preceq_{\epsilon+\delta+\epsilon\delta} X'\times W'$.

$iv)\ \ X\preceq_\epsilon X'\wedge Y\preceq_{\delta}Y'\ and\ X',Y'\ are\ separable\ by\ a\ circuit,\ then\ X\cup Y\preceq_{\epsilon+\delta}X'\cup Y'$.

\end{proposition}

\begin{proposition}[Je\v{r}\'{a}bek \cite{Japx}](in \APC)\label{really} 

\noindent 1.\ \ Let $X,Y\subseteq 2^n$ be definable by circuits, $s,t,u\leq 2^n$, $\epsilon, \delta,\theta, \gamma < 1, \gamma^{-1}\in Log $. Then

i)\ \ $X\preceq_{\gamma} Y$ or $Y\preceq_{\gamma} X$,

ii)\ \ $s\preceq_\epsilon X\preceq_\delta t\Rightarrow s<t+(\epsilon+\delta+\gamma)2^n$,

iii)\ \ $X\preceq_{\epsilon} Y\Rightarrow 2^n\backslash Y\preceq_{\epsilon+\gamma} 2^n\backslash X $,

iv)\ \ $X\approx_{\epsilon} s\wedge Y\approx_{\delta} t\wedge X\cap Y\approx_{\theta} u\Rightarrow X\cup Y\approx_{\epsilon+\delta+\theta+\gamma} s+t-u$.

\medskip

\noindent 2. (Disjoint union) Let $X_i\subseteq 2^n$, $i<m$ be defined by a sequence of circuits and $\epsilon,\delta\leq 1$, $\delta^{-1}\in Log$. If $X_i\preceq_\epsilon s_i$ for every $i<m$, then $\bigcup_{i<m} (X_i\times \{i\})\preceq_{\epsilon+\delta} \sum_{i<m} s_i$.



\def\chernoff{\noindent 3. (Chernoff's bound) Let $X\subseteq 2^n, m\in Log, 0\leq \epsilon,\delta,p\leq 1$ and $X\succeq_{\epsilon} p2^n$. Then $$\{w\in (2^n)^m|\ |\{i<m|w_i\in X\}|\leq m(p-\delta)\}\preceq_0 c4^{m(c\epsilon-\delta^2)}2^{nm}$$ for some constant $c$, where $w$ is treated as a sequence of $m$ numbers less than $2^n$ and $w_i$ is its $i$-th member.}
\end{proposition}





When proving $\Sigma^b_1$ statements in \APC we can afford to work in $\SB+dWPHP(\pv)+BB(\Sigma^b_2)$ and, in fact, assuming the existence of a single hard function in \PV gives us the full power of \APC. Here, $BB(\Sigma^b_2)$ is defined as $BB(\Sigma^b_1)$ but with $A\in\Sigma^b_2$. 

\begin{lemma}[\cite{MP}]\label{singlehard} Suppose $\SB+dWPHP(\pv)+BB(\Sigma^b_2)\vdash \exists y A(x,y)$ for $A\in\Sigma^b_1$. Then, for every $\epsilon<1$, there is $k$ and $\pv$-functions $g,h$ such that \PV proves $$|f|\ge |x|^k\wedge \exists y, |y|=||f||, C_h(y)\ne f(y)\rightarrow A(x,g(x,f))$$ where $f(y)$ is the $y$th bit of $f$, $f(y)=0$ for $y>|f|$, and $C_h$ is a circuit of size $\leq 2^{\epsilon ||f||}$ generated by $h$ on $f,x$. Moreover, $\APC\vdash\exists y A(x,y)$.
\end{lemma}

\def\prfofMPlemma{
\proof By \cite[Corollary 4.12]{Jdual}, $\SB+dWPHP(\pv)+BB(\Sigma^b_2)\vdash\exists y A(x,y)$ implies $\SB+dWPHP(\pv)\vdash\exists y A(x,y)$. Then, following Thapen's proof of \cite[Theorem 4.2]{Tstr} (cf. also \cite[Proposition 1.14]{Jdual}), there is $\ell$ and $h\in\pv$ such that \SB proves $$(\forall v\leq 2^{8|x|^\ell}\exists u\leq 2^{4|x|^\ell},\ h(u)=v)\vee \exists y A(x,y).$$ By Buss's witnessing theorem it now suffices to show that for every $\epsilon<1$ there is $k$ such that \SB proves \begin{multline*}(\forall v\leq 2^{8|x|^\ell}\exists u\leq 2^{4|x|^\ell},\ h(u)=v)\rightarrow\\ 
(|f|\ge |x|^k\rightarrow \exists\text{ circuit }C\text{ of size} \leq 2^{\epsilon||f||}\ \forall y, |y|=||f||, C(y)=f(y)).\end{multline*} Argue in \SB. The surjection $h:2^m\rightarrow 2^{2m}$, where $m=4|x|^\ell$, is computed by a circuit of size $m^{\ell'}$ for a standard $\ell'$. Following Je\v{r}\'{a}bek's \SB-proof of \cite[Proposition  3.5]{Jdual}, this implies that every (number) $f$ viewed as a truth-table of length $|f|$ is computed by a size $O(m|m|+m^{\ell'}|\lceil |f|/m \rceil|)$ circuit with $||f||$ inputs. 
For sufficiently large $k$, $|f|\ge |x|^k$ implies that this size is $\le 2^{\epsilon||f||}$.

The ``moreover'' part is a consequence of $\APC\vdash \forall n\in LogLog\ \exists f:2^n\rightarrow 2, \LBtt(f,2^{n/4})$, cf. \cite[Corollary 3.3]{Jdual}. 
 \qed}

\def\inverselemma{
\bigskip

Lemma \ref{singlehard} allows us to use the $BB(\Sigma^b_2)$ collection scheme for proving $\Sigma^b_1$-statements in \APC. Unfortunately, when collecting circuits witnessing $\preceq_\epsilon$ predicates given by $\exists\Pi^b_2$ formulas the $BB(\Sigma^b_2)$ collection is a priori not sufficient. To overcome this complication the quantifier complexity of $\preceq_\epsilon$ can be pushed down to $\Sigma^b_2$ because the circuits counting sizes of sets in \APC are invertible.

\begin{lemma}[\cite{MP}]\label{invert}(in \APC) Let $X\subseteq 2^n$ be defined by a circuit and $\epsilon^{-1}\in Log$. Suppose $X\preceq_\epsilon s$. Then, $X\preceq_\epsilon s+3\epsilon2^n$ is expressible by a provable $\Sigma^b_2$ formula.
\end{lemma}

}

\def\pvine{
\subsection{Standard inequalities in \PV}\label{spv}

For a \pv-function symbol $f$ and $n\in Log$, in \PV we can define inductively $\sum_{i=0}^n f(i)$. Similarly, we can define iterated products, factorials, and binomial coefficients. It is easy to see that, by induction, \PV proves: $n\in Log\rightarrow \sum_{i=0}^n {n\choose i}=2^n$. 

\begin{proposition}[Stirling's bound, cf. Je\v{r}\'{a}bek \cite{Jdual}]\label{stirling} There is a $c>1$ such that \PV proves: $$0<k<n\in Log\rightarrow \frac{1}{c}{n \choose k}<\frac{n^n}{k^k(n-k)^{n-k}}\left(\left\lfloor\sqrt{\frac{k(n-k)}{n}}\right\rfloor+1\right)^{-1}<c {n \choose k}.$$
\end{proposition}

\begin{proposition}\label{concentration} For each $\epsilon>0$ there is an $n_0$ such that \PV proves: $$n_0<n\in Log\rightarrow \sum_{i=0}^{\lfloor n/2+n^{1/3}\rfloor} {n\choose i}<\left(\frac{1}{2}+\epsilon\right)2^n.$$
\end{proposition}

\proof $\sum_{i=0}^{\lfloor n/2\rfloor-1} {n\choose i}=\frac{1}{2}\left(\sum_{i=0}^{\lfloor n/2\rfloor-1} {n\choose i}+\sum_{i=0}^{\lfloor n/2\rfloor-1} {n\choose n-i}\right)<2^{n-1}$ and by Stirling's bound, for some constant $c>1$, $$\sum_{i=\lfloor n/2\rfloor}^{\lfloor n/2+n^{1/3}\rfloor} {n\choose i}<(n^{1/3}+1){n\choose \lfloor n/2\rfloor}<2^n4c\left(\frac{n^{1/3}}{\lfloor n^{1/2}/2\rfloor}+\frac{1}{\lfloor n^{1/2}/2\rfloor}\right)$$ where to verify the last inequality for odd $n$ we used also the provability of $a,b\in Log$, $b>0\rightarrow (1+a/b)\leq 4^{a/b}$ shown in \cite[Stirling's bound, Claim 1]{Jdual}. \qed

\begin{proposition}\label{expbound} \PV proves: $$a,b\in Log,\ b>a+1\rightarrow (b-a)^b\leq b^b/2^a.$$ Note that the conclusion implies $(1-a/b)\leq 2^{-a/b}$.
\end{proposition}

\proof Proceed as in the proof of Claim 2 in the proof of Stirling's bound \cite{Jdual} but instead of Claim 1 use the inequality $b^b\leq  (b+1)^b/2$. \qed}

\subsection{Formalizing complexity-theoretic statements}\label{s:formalizations}

\subsubsection{Circuit lower bounds} An `almost everywhere' circuit lower bound for circuits of size $s$ and a function $f$ says that for every sufficiently big $n$, for each circuit $C$ with $n$ inputs and size $s$, there exists an input $y$ on which the circuit $C$ fails to compute $f(y)$.

If $f:\{0,1\}^n\rightarrow \{0,1\}$ is an \NPtime function and $s=n^k$ for a constant $k$, this can be written down as a $\forall\Sigma^b_2$ formula $\LB(f,n^k)$, $$\forall n,\ n>n_0\ \forall\ \text{circuit}\ C\ \text{of size}\ \leq n^k\ \exists y,\ |y|=n,\ C(y)\neq f(y),$$ where $n_0$ is a constant and $C(y)\neq f(y)$ is a $\Sigma^b_2$ formula stating that a circuit $C$ on input $y$ outputs the opposite value of $f(y)$. The intended meaning of `$\exists y, |y|=n$' is to say that $y$ is a string from $\{0,1\}^n$. This is a slight abuse of notation as, formally, $|y|=n$ fixes the first bit of $y$ to 1.

If we want to express $s(n)$-size lower bounds for $s(n)$ as big as $2^{O(n)}$, we add an extra assumption on $n$ stating that $\exists x,\ n=||x||$. That is, the resulting formula $\LBtt(f,s(n))$ has form `$\forall x,n; n=||x||\rightarrow\dots$'. Treating $x,n$ as free variables, $\LBtt(f,s(n))$ is $\Pi^b_1$ if $f$ is, for instance, \SAT because $n=||x||$ implies that the quantifiers bounded by $2^{O(n)}$ 
are sharply bounded. Moreover, allowing $f\in \mathsf{NE}$ lifts the complexity of $\LBtt(f,s(n))$ just to $\forall\Sigma^b_1$. The function $s(n)$ in $\LBtt(f,s(n))$ is assumed to be a \pv-function with input $x$ (satisfying $||x||=n$).
\medskip

In terms of the $Log$-notation, $\LB(f,n^k)$ implicitly assumes $n\in Log$ while $\LBtt(f,n^k)$ assumes $n\in LogLog$.
By chosing the scale of $n$ we are determining how big objects are going to be `feasible' for theories reasoning about the statement. In the case $n\in LogLog$, the truth-table of $f$ (and everything polynomial in it) is feasible. Assuming just $n\in Log$ means that only the objects of polynomial-size in the size of the circuit are feasible. Likewise, the theory reasoning about the circuit lower bound is less powerful when working with $\LB(f,n^k)$ than with $\LBtt(f,n^k)$. (The scaling in $\LBtt(f,s)$ corresponds to the choice of parameters in natural proofs and in the formalizations by Razborov \cite{Rup}.)
\medskip

We can analogously define formulas $\LBtt(f,s(n),t(n))$ expressing an average-case lower bound for $f$, where $f$ is a free variable (with $f(y)$ being the $y$th bit of $f$ and $f(y)=0$ for $y>|f|$). More precisely, $\LBtt(f,s(n),t(n))$ generalizes $\LBtt(f,s(n))$ by saying that each circuit of size $s(n)$ fails to compute $f$ on at least $t(n)$ inputs, for \pv-functions $s(n),t(n)$. Since $n\in LogLog$, $\LBtt(f,s(n),t(n))$ is $\Pi^b_1$.
\medskip



\noindent {\bf Propositional version.} An $s(n)$-size circuit lower bound for a function $f:\{0,1\}^n\rightarrow \{0,1\}$ can be expressed by a $poly(2^n,s)$-size propositional formula ${\sf tt}(f,s)$, $$\bigvee_{y\in\{0,1\}^n} f(y)\neq C(y)$$ where the formula $f(y)\neq C(y)$ says that an $s(n)$-size circuit $C$ represented by $poly(s)$ variables does not output $f(y)$ on input $y$. The values $f(y)$ are fixed bits. That is, the whole truth-table of $f$ is hard-wired in ${\sf tt}(f,s)$. 

The details of the encoding of the formula $\ttable(f,s)$ are not important for us as long as the encoding is natural because systems like \EF considered in this paper can reason efficiently about them. We will assume that $\ttable(f,s)$ is the formula resulting from the translation of $\Pi^b_1$ formula $\LBtt(h,s)$, where $n_0=0$, $n,x$ are substituted after the translation by fixed constants so that $x=2^{2^n}$, and $h$ is a free variable (with $h(y)$ being the $y$th bit of $h$ and $h(y)=0$ for $y>|h|$) which is substituted after the translation by constants defining $f$.

\def\erasepartfor{
Even more feasible, uniform, encoding follows from translations of $\LB(f,n^k)$. This requires an efficient witnessing of existential quantifiers in $\LB(f,n^k)$ collapsing its complexity to $\forall\Sigma^b_0$. Such a p-time witnessing of $\LB(\SAT,n^k)$ follows, for example, from the existence of one-way permutations and a function in $\mathsf{E}$ hard for subexponential-size circuits, cf. \cite[Proposition 4.3]{clba} \footnote{Proposition 4.3 in \cite{clba} shows just the existence of an S-T protocol witnessing $\LB(\SAT,n^k)$ but the p-time witnessing easily follows. 
}. Further, by the KPT theorem \cite{KPT}, whenever $\PV\vdash\LB(f,n^k)$ we get a sequence of finitely many p-time functions $\overline{w}=w_1,\dots,w_c$ witnessing the existential quantifiers in $\LB(f,n^k)$. 
$\LB(f,n^k)$ witnessed by $\overline{w}$ can be equivalently expressed by a sequence of $poly(n)$-size propositional formulas ${\sf lb}_{\overline{w}}(f,n^k)$. 
}

Analogously, we can express average-case lower bounds by propositional formulas $\ttable(f,s(n),t(n))$ obtained by translating $\LBtt(h,s(n),t(n)2^n)$, with $n_0=0$, fixed $x=2^{2^n}$ and $h$ substituted after the translation by $f$.

\subsubsection{Learning algorithms}\label{ss:lform}
A circuit class $\mathcal{C}$ is defined by a \pv-formula if there is a \pv-formula defining the predicate $C\in\mathcal{C}$. Definition \ref{d:lear} can be formulated in the language of ${\sf HARD}^A$: A circuit class $\mathcal{C}$ (defined by a \pv-formula) is learnable over the uniform disribution by a circuit class $\mathcal{D}$ (defined by a \pv-formula) up to error $\epsilon$ with confidence $\delta$, 
if there are randomized oracle circuits $L^f$ from $\mathcal{D}$ such that for every Boolean function $f:\{0,1\}^n\mapsto\{0,1\}$ (represented by its truth-table) computable by a circuit from $\mathcal{C}$, for each $\gamma^{-1}\in Log$, when given oracle access to $f$, input $1^n$ and the internal randomness $w \in \{0,1\}^*$, $L^f$ outputs the description of a circuit satisfying 
  \begin{equation*}
    \Pr_w [ L^f(1^n, w) \text{ } (1-\epsilon) \text{-approximates } f ]_{\gamma} \geq \delta.
  \end{equation*}
The inner probability of approximability of $f$ by $L^f(1^n,w)$ is counted exactly. This is possible because $f$ is represented by its truth-table, which implies that $2^n\in Log$.

\medskip
\noindent {\bf Propositional version.} 
In order to translate the definition of learning algorithms to propositional formulas and to the language of \PV we need to look more closely at the definition of  ${\sf HARD}^A$.


\PV can be relativized to $\PV(\alpha)$. The new function symbol $\alpha$ is then allowed in the inductive clauses for introduction of new function symbols. 
This means that the language of $\PV(\alpha)$, denoted also $\pv(\alpha)$, contains symbols for all p-time oracle algorithms.

\begin{proposition}[Je\v{r}\'{a}bek \cite{Jwphp}]\label{p:hardf} For every constant $\epsilon <1/3$ there exists a constant $n_0$ such that \APC proves: for every $n\in LogLog$ such that $n>n_0$, there exist a function $f:2^n\rightarrow 2$ such that no circuit of size $2^{\epsilon n}$ computes $f$ on $\ge (1/2+1/2^{\epsilon n})2^n$ inputs. \end{proposition}

\begin{definition}[Je\v{r}\'{a}bek \cite{Jwphp}]\label{d:hard} The theory ${\sf HARD}^A$ is an extension of the theory $\PV(\alpha)+dWPHP(\pv(\alpha))$ by the axioms
\smallskip

1. $\alpha(x)$ is a truth-table of a Boolean function in $||x||$ variables,
 
2. $\LBtt(\alpha(x),2^{||x||/4},2^{||x||}(1/2-1/2^{||x||/4}))$, for constant $n_0$ from Proposition \ref{p:hardf},

3. $||x||=||y||\rightarrow \alpha(x)=\alpha(y)$.
\end{definition}

By inspecting the proof of Lemma 2.14 in \cite{Japx}, we can observe that on each input $C,2^n,2^{\epsilon^{-1}}$ the $\PV(\alpha)$-function $Size$ calls $\alpha$ just once (to get the truth-table of a hard function  which is then used as the base function of the Nisan-Wgiderson generator). In fact, $Size$ calls $\alpha$ on input $x$ which  depends only on $|C|$, the number of inputs of $C$ and w.l.o.g. also just on $|\epsilon^{-1}|$ (since decreasing $\epsilon$ leads only to a better approximation). In combination with the fact that the approximation $Size(C,2^n,2^{\epsilon^{-1}})\approx_{\epsilon} X$, for $X\subseteq 2^n$ defined by $C$, is not affected by a particular choice of the hard boolean function generated by $\alpha$, we get that \APC  proves $$\LBtt(y,2^{||y||/4},2^{||y||}(1/2-1/2^{||y||/4}))\wedge ||y||=S(C,2^n,2^{\epsilon^{-1}}) \rightarrow\ Sz(C,2^n,2^{\epsilon^{-1}},y)\approx_{\epsilon} X,$$ where $Sz$ is defined as $Size$ with the only difference that the call to $\alpha(x)$ on $C,2^n,2^{\epsilon^{-1}}$ is replaced by $y$ and $S(C,2^n,2^{\epsilon^{-1}})=||x||$ for a \pv-function $S$. ($S$ is given by a subcomputation of $Size$ specifying $||x||$, for $x$ on which $Size$ queries $\alpha(x)$.)

This allows us to express $\Pr_{x<t}[x\in X]_{\epsilon}=a$, where $\epsilon^{-1}\in Log$ and $X\cap t\subseteq 2^{|t|}$ is defined by a circuit $C$, without a $\PV(\alpha)$ function, by formula
 $$\forall y\ (\LBtt(y,2^{||y||/4},2^{||y||}(1/2-1/2^{||y||/4}))\wedge ||y||=S(C,2^{|t|},2^{\epsilon^{-1}})\rightarrow Sz(C,2^{|t|},2^{\epsilon^{-1}},y)/t=a).$$ 
We denote the resulting formula by $\Pr^y_{x<t}[x\in X]_{\epsilon}=a$. We will use the notation $\Pr^y_{x<t}[x\in X]_{\epsilon}$ in equations with the intended meaning that the equation holds for the value $Sz(\cdot,\cdot,\cdot,\cdot)/t$ under corresponding assumptions. For example, $t\cdot \Pr^y_{x<t}[x\in X]_{\epsilon}\preceq_{\delta} a$ stands for `$\forall y, \exists v,\exists$ circuit $\hat{C}$ (defining a surjection) which witnesses that $\LBtt(y,2^{||y||/4},2^{||y||}(1/2-1/2^{||y||/4}))\wedge ||y||=S(C,2^{|t|},2^{\epsilon^{-1}})$ implies $Sz(C,2^{|t|},2^{\epsilon^{-1}},y)\preceq_{\delta} a$'.

The definition of learning can be now expressed without a $\PV(\alpha)$ function: If circuit class $\mathcal{C}$ is defined by a \pv-function, the statement that a given oracle algorithm $L$ (given by a \pv-function with oracle queries) learns a circuit class $\mathcal{C}$ over the uniform distribution up to error $\epsilon$ with confidence $\delta$ can be expressed as before with the only difference that we replace 
    $\Pr_w [ L^f(1^n, w) \text{ } (1-\epsilon) \text{-approximates } f ]_{\gamma} \geq \delta$ 
by 
 \begin{equation*}
    \Pr^y_w [ L^f(1^n, w) \text{ } (1-\epsilon) \text{-approximates } f ]_{\gamma} \geq \delta. 
  \end{equation*}

Since the resulting formula $A$ defining learning  is not $\Pi^b_1$ (because of the assumption $\LBtt$) 
 we cannot translate it to propositional logic. We will sidestep the issue by translating only the formula $B$ obtained from $A$ by deleting subformula $\LBtt$ (but leaving $||y||=S(\cdot,\cdot,\cdot)$ intact) and replacing the variables $y$ by  fixed bits representing a hard boolean function. In more detail, $\Pi^b_1$ formula $B$ can be translated into a sequence of propositional formulas $\mathsf{lear}^y_{\gamma}(L,\mathcal{C},\epsilon,\delta)$ expressing that ``if $C\in\mathcal{C}$ is a circuit computing $f$, then $L$ querying $f$ generates a circuit $D$ such that $\Pr[D(x)=f(x)]\ge 1-\epsilon$ with probability $\ge \delta$, which is counted approximately with precision $\gamma$". Note that $C,f$ are represented by free variables and that there are also free variables for error $\gamma$ from approximate counting and for 
 Boolean functions $y$. 
As in the case of \ttable-formulas, we fix $|f|=2^n$, so $n$ is not a free variable. %
Importantly, $\mathsf{lear}^y_{\gamma}(L,\mathcal{C},\epsilon,\delta)$ does not postulate that $y$ is a truth-table of a hard boolean function.
 Nevertheless, for any fixed (possibly non-uniform) bits representing a sequence of Boolean functions
 $h=\{h_m\}_{m>n_0}$
 such that $h_m$ is not $(1/2+1/2^{m/4})$-approximable by any circuit of size $2^{m/4}$, we can obtain formulas $\mathsf{lear}^h_{\gamma}(L,\mathcal{C},\epsilon,\delta)$ by substituting bits $h$ for $y$. 

Using a single function $h$ in $\mathsf{lear}^h_{\gamma}(L,\mathcal{C},\epsilon,\delta)$ does not ruin the fact that (the translation of function) $Sz$ approximates the respective probability with accuracy $\gamma$ because $Sz$ queries a boolean function $y$ which depends just on the number of atoms representing $\gamma^{-1}$ 
 and on the size of the circuit $D$ defining the predicate we count together with the number of inputs of $D$. The size of $D$ and the number of its inputs are w.l.o.g. determined by the number of inputs of $f$. 

If we are working with formulas $\mathsf{lear}^h_{\gamma}(L,\mathcal{C},\epsilon,\delta)$, where $h$ is a sequence of bits representing a hard boolean function, in a proof system which cannot prove efficiently that $h$ is hard, our proof system might not be able to show that the definition is well-behaved - it might not be able to derive some standard properties of the function $Sz$ 
 used inside the formula. Nevertheless, in our theorems this will never be the case: our proof systems will always know that $h$ is hard.

In formulas $\mathsf{lear}^y_{\gamma}(L,\mathcal{C},\epsilon,\delta)$ we can allow $L$ to be a sequence of nonuniform circuits, with a different advice string for each input length. One way to see that is to use additional input to $L$ in $\Pi^b_1$ formula $B$, then translate the formula to propositional logic and substitute the right bits of advice for the additional input. Again, the precise encoding of the formula $\mathsf{lear}^y_{\gamma}(L,\mathcal{C},\epsilon,\delta)$ does not matter very much to us but in order to simplify proofs we will assume that $\mathsf{lear}^y_{\gamma}(L,\Circuit[n^k],\epsilon,\delta)$ has the from $\neg\ttable(f,n^k)\rightarrow R$, where $n, k$ are fixed, 
$f$ is represented by free variables and $R$ is the remaining part of the formula expressing that $L$ generates a suitable circuit with high probability. 

\subsubsection{Nisan-Wigderson generators}

The core theorem underlying approximate counting in \APC is the following formalization of Nisan-Wigderson generators (NW), cf. \cite[Proposition 4.7]{Jwphp}.

\begin{theorem}[Je\v{r}\'{a}bek \cite{Jwphp}]\label{t:nwpv} Let $0<\gamma<1$. There are constants $c>1$ and $\delta'>0$ such that for each $\delta<\delta'$ there is a $poly(2^m)$-time function $$NW:\{0,1\}^{cm}\times\{0,1\}^{2^m}\mapsto\{0,1\}^{\lfloor 2^{\delta m}\rfloor}$$ such that \SB proves: ``If $2^m\in Log$ and $f:\{0,1\}^{m}\mapsto \{0,1\}$ is a Boolean function such that no circuit of size $2^{\epsilon m}$ computes $f$ on $> (1/2+1/2^{\epsilon m})2^m$ inputs, then for each $(2^{\epsilon m}-\lceil 2^{(\delta+\gamma)m}\rceil-1)$-size circuit $D$ with $\lfloor 2^{\delta m}\rfloor$ inputs, $$2^{\lfloor 2^{\delta m}\rfloor}\times \{z<2^{cm}\mid D(NW_f(z))=1\}\succeq_{e} 2^{cm}\times \{x<2^{\lfloor 2^{\delta m}\rfloor}\mid D(x)=1\},$$ where $e:=\lceil 2^{\delta m}\rceil/2^{\epsilon m}$ and
$NW_f(z):=NW(z,f)$." \end{theorem}

Theorem \ref{t:nwpv} shows that $\Pr_x^y[D(x)=1]_{\theta}$ is \SB-provably similar to $\Pr_z^y[D(NW_f(z))=1]_{\theta}$,  for $\theta^{-1}\in Log$. To see this, note that $$2^{cm}\Pr_z^y[D(NW_f(z))=1]_{\theta}\approx_{\theta}  \{z<2^{cm}\mid D(NW_f(z))=1\}.$$ Hence, by Proposition \ref{lem} $iii)$,
$$2^{\lfloor 2^{\delta m}\rfloor }2^{cm}\Pr_z^y[D(NW_f(z))=1]_{\theta}\succeq_{\theta} 2^{\lfloor 2^{\delta m}\rfloor} \times \{z<2^{cm}\mid D(NW_f(z))=1\}.$$ Similarly, $$2^{cm}\times \{x<2^{\lfloor 2^{\delta m}\rfloor}\mid D(x)=1\}\succeq_{\theta}2^{\lfloor 2^{\delta m}\rfloor }2^{cm}\Pr_x^y[D(x)=1]_{\theta}.$$

Therefore, by Proposition \ref{lem} $ii)$, the conclusion of Theorem \ref{t:nwpv} implies
$$2^{\lfloor 2^{\delta m}\rfloor}2^{cm}\Pr_z^y[D(NW_f(z))=1]_{\theta}\succeq_{2\theta+e} 2^{\lfloor 2^{\delta m}\rfloor}2^{cm}\Pr_x^y[D(x)=1]_{\theta}.$$
If the size of $D$ is $\le 2^{\epsilon m}-\lceil 2^{(\delta+\gamma)m}\rceil-2$, the same inequality holds for $\neg D$ instead of $D$.

\section{Feasible anticheckers.}\label{s:fant}

If there is an $n^{k}$-size circuit computing $\SAT_n$, there is a $poly(n^{k})$-size circuit $B$ with $n$ inputs and $\le n$ outputs such that $\forall x,y\in \{0,1\}^n, (\SAT_n(x,y)\rightarrow \SAT_n(x,B(x)))$. We use this to formulate the existence of anticheckers for \SAT as a $\forall\Pi^b_1$ statement.

\begin{theorem}[`CC $\leftarrow$ PC' from feasible anticheckers]\label{t:fant}\ \newline
Let $k\ge 3$ be a constant and assume that there is a p-time function $f$ such that \SB proves:
\begin{itemize}
\item[] ``$\forall 1^n$, $f(1^n)$ is a $poly(n^{k})$-size circuit $B$ such that $$\forall x,y \in \{0,1\}^n, [\SAT_n(x,y)\rightarrow \SAT_n(x,B(x))]$$ or $\big(f(1^n)$ outputs sets $A^{\SAT_n,n^k}_n,A'\subseteq \{0,1\}^n$, $D\subseteq A^{\SAT_n,n^k}_n\times A'$ of size $poly(n^k)$ such that 
$$\forall x\in A^{\SAT_n,n^k}_n[\exists y_x\in A', \left<x,y_x\right>\in D\wedge \forall z,y\in A', (\left<x,y\right>\in D\wedge \left<x,z\right>\in D\rightarrow y=z)]$$ and
$\forall\ n^{k}$-size circuit $C$, $$\forall x\in A^{\SAT_n,n^k}_n \forall y\in\{0,1\}^n\ [\SAT_n(x,y)\rightarrow \SAT_n(x,y_x)]\wedge$$ $$\exists x\in A^{\SAT_n,n^k}_n, \SAT_n(x,y_x)\ne C(x)\big)."$$
\end{itemize}
Then, proving that \EF is not p-bounded implies $\SAT_n\notin\Circuit[n^{k}]$ for infinitely many $n$.
\end{theorem}

\proof
The statement assumed to have an \SB-proof is $\forall\Pi^b_1$, 
so there are p-size \EF-proofs of its propositional translation. If we now assume that $\exists n_0\forall n>n_0, \SAT_n\in \Circuit [n^k]$, there are circuits $C$ and $y\in \{0,1\}^n$ falsifying the second disjunct of the translated assumption for $n>n_0$. Consequently, \EF proves efficiently that the circuits generated by $f(1^n)$ solve $\SAT_n$, which implies that \EF is p-bounded.
\qed

\medskip

\medskip

\noindent {\bf Existential quantifiers instead of witnessing.} If we used the existential quantifiers instead of function $f$ in Theorem \ref{t:fant}, the resulting statement $S$ formalizing the existence of anticheckers would be $\forall\Sigma^b_2$. By the KPT theorem \cite{KPT}, \PV-provability of $S$ would then imply the existence of p-time functions $f_1,\dots,f_{c}$, for a constant $c$, with a \PV-proof of:
\begin{itemize}
\item[] ``$\forall 1^n$, $\forall x^1,\dots,x^c,y^1,\dots,y^c, \tilde{y}^1,\dots,\tilde{y}^c \in \{0,1\}^n, \forall n^k$-size circuits $C^1,\dots, C^c$,\\ $f_1(1^n)$ outputs a $poly(n^{k})$-size circuit $B$ and $A^{\SAT_n,n^k}_n,A'\subseteq \{0,1\}^n, D\subseteq A^{\SAT_n,n^k}_n\times A'$ of size $poly(n^k)$ such that the following predicate, denoted $P_{f_1}(x^1,y^1,C^1,\tilde{y}^1)$, holds: $$
\big(\SAT_n(x^1,y^1)\rightarrow \SAT_n(x^1,B(x^1))\big)\vee \big(D'(\tilde{y}^1)\wedge \exists \tilde{x}\in A^{\SAT_n,n^k}_n, \SAT_n(\tilde{x},y_{\tilde{x}})\ne C(\tilde{x})\big),$$ where $D'(\tilde{y}^1)$ stands for the remaining part of the $\Sigma^b_0$ subformula of $S$,\\ or $f_2(1^n, x^1,y^1,C^1,\tilde{y}^1)$ outputs a $poly(n^{k})$-size circuit $B$ and $A^{\SAT_n,n^k}_n,A',D$ of size $poly(n^k)$ such that $P_{f_2}(x^2,y^2,C^2,\tilde{y}^2)$ holds, or \\
$\dots$\\
or $f_c(1^n,x^1,\dots,x^{c-1},y^1,\dots,y^{c-1},C^1,\dots,C^{c-1},\tilde{y}^1,\dots,\tilde{y}^{c-1})$ outputs a $poly(n^{k})$-size circuit $B$ and $A^{\SAT_n,n^k}_n,A', D$ of size $poly(n^k)$ such that $P_{f_c}(x^c,y^c,C^c,\tilde{y}^c)$."
\end{itemize}
The resulting $\forall\Pi^b_1$-statement could be translated to propositional tautologies with p-size \EF-proofs. 
However, given $\forall n, \SAT_n\in \Circuit [n^k]$, we could not directly obtain p-size \EF-proofs of tautologies stating that one of the functions $f_1,\dots,f_c$ generates a circuit solving $\SAT_n$. This is because $B$ and $A^{\SAT_n,n^k}_n$ generated by $f_2$ depend on $y^1$. For the same reason, it seems possible for \EF to prove efficiently $\SAT\in \Ppoly$ (using the formalization based on the KPT witnessing) without proving efficiently all tautologies. 

\section{One-way functions from  $\NP\not\subseteq\Ppoly$}\label{s:owpnp}

\begin{theorem}[`CC $\leftarrow$ PC' from `OWF $\leftarrow$ $\NP\not\subseteq\Ppoly$' \& hardness of {\sf E}]\label{t:owpnp}\ \newline
Assume that for each sufficiently big $n$, each $2^{n/4}$-size circuit fails to compute $h'\in {\sf E}$ on $\ge (1/2-1/2^{n/4})$ of all inputs. Further, assume that there is a p-time function $h:\{0,1\}^n\mapsto\{0,1\}^{u(n)}$ such that for each constants $c,d$, there is a p-time function $f_2$ and a constant $0<\epsilon<1$ such that \SB proves:
\begin{itemize}
\item[]  ``$\forall n, \forall\ cn^{c}$-size circuit $C$ with $u(m)$ inputs and $m$ outputs such that $n\le dm^{d}$,\\
$\big( f_2(C)$ is a $poly(n)$-size circuit $B$ such that $$\forall x,y \in \{0,1\}^{\lfloor n^{\epsilon}\rfloor}, [\SAT_{\lfloor n^{\epsilon}\rfloor}(x,y)\rightarrow \SAT_{\lfloor n^{\epsilon}\rfloor}(x,B(x))]$$ or  $$\Pr^{y}_{x\in \{0,1\}^{m}}[h(C(h(x)))=h(x)]_{\frac{1}{m}}<1/2\big)."$$
\end{itemize}
Then, proving that $\EF+\ttable(h',2^{n/4},1/2-1/2^{n/4})$ is not p-bounded implies $\SAT\notin\Ppoly$.
\end{theorem}

The system $\EF+\ttable(h',2^{n/4},1/2-1/2^{n/4})$ is defined in the same way as in the introduction. That is, $\EF+\ttable(h',2^{n/4},1/2-1/2^{n/4})$ is an extension of \EF which can use substitutional instances of $\ttable(h'_n,2^{n/4},1/2-1/2^{n/4})$, for sufficiently big $n$, in its proofs.
\bigskip

The proof of Theorem \ref{t:owpnp} is based on the following lemma formalizing a conditional witnessing of $\NP\not\subseteq \Ppoly$, cf. \cite{Pclba,MP}.

\begin{lemma}[Formalized witnessing of $\NP\not\subseteq \Ppoly$ from OWF \& hardness of {\sf E}]\label{l:owpwit}\ \newline
Let $k\ge 1$ be a constant. For each p-time functions $h:\{0,1\}^n\mapsto \{0,1\}^{u(n)}$ and $f_1$, there are p-time functions $f_0,f_{-1},f_{-2}$ and constants $b, n_1$ such that $\SB+dWPHP(\pv)$ proves: ``$\forall 1^n>n_1$, $\forall m$ such that $n/2^b\le 2^{bm}\le n$, if $$\LBtt'(f_1(1^{2^{m}}),2^{m/4},2^{m}(1/2-1/2^{m/4})),$$

then $f_0(1^n,m)$ outputs sets $A,A'\subseteq \{0,1\}^n$ of size $poly(n)$ such that $$\forall x\in A\ \exists y_x\in A'\ \SAT_n(x,y_x)$$ and $\big(\forall\ n^{k}$-size circuit $C$ with $n$ inputs and $\le n$ outputs, $$\exists x\in A, \neg\SAT_n(x,C(x))$$ or $f_{-1}(C,m)$ outputs a $poly(n)$-size circuit $C'$ with $u(n')$ inputs and $n'$ outputs, where $n\le f_{-2}(1^{n'})$, such that $$\Pr^{y}_{x\in \{0,1\}^{n'}}[h(C'(h(x)))=h(x)]_{\frac{1}{n'}}\ge 1/2\big)."$$ Here, $\LBtt'$ is obtained from $\LBtt$ by setting $m_0=0$ and skipping the universal quantifier on $m$ (so $m$ in $\LBtt'$ is the same as the universally quantified $m$ in the \SB-provable statement).
\end{lemma}

\proof[Proof of Theorem \ref{t:owpnp} from Lemma \ref{l:owpwit}] Intuitively, Theorem \ref{t:owpnp} assumes that the hardness of \SAT yields a function $h$ which is hard to invert. Lemma \ref{l:owpwit} shows that such $h$ can be used to find an error of each small circuit attempting to compute \SAT. Combining the assumption of Theorem \ref{t:owpnp} with Lemma \ref{l:owpwit} we obtain a p-time function $f$ such that for each small circuit $C$, either $C$ solves \SAT or $f(C)$ finds an error of $C$. Moreover, this holds provably in $\SB+dWPHP(\pv)$, so the propositional translation of the correctness of the witnessing statement has short proofs in $\EF+\ttable(h',2^{n/4},1/2-1/2^{n/4})$. This will allow us to derive the desired implication similarly as in the proof of Theorem \ref{t:core}.

We proceed with a formal proof.
\bigskip

The assumption of Theorem \ref{t:owpnp} in combination with Lemma \ref{l:owpwit} implies that for each $k\ge 1$, for p-time $f_1$ generating the truth-table of $h'$, there is $0<\epsilon<1$ and $b, n_1$ such that $\SB+dWPHP(\pv)$ proves the following statement $S$:
\medskip

\noindent ``$\forall 1^n>n_1$, $\forall m, n/2^b\le 2^{bm}\le n$, if $$\LBtt'(f_1(1^{2^{m}}),2^{m/4},2^{m}(1/2-1/2^{m/4})),$$ then $f_0(1^n,m)$ outputs $A,A'\subseteq \{0,1\}^n$ such that $$\forall x\in A\ \exists y_x\in A'\ \SAT_n(x,y_x)$$ and $\big(\forall n^{k}$-size circuit $C$ with $n$ inputs and $\le n$ outputs, $$\exists x\in A, \neg\SAT_n(x,C(x))$$ or
$f_2(f_{-1}(C,m))$ outputs a circuit $B$ such that $$\forall x,y \in \{0,1\}^{\lfloor n^{\epsilon}\rfloor}, [\SAT_{\lfloor n^{\epsilon}\rfloor}(x,y)\rightarrow \SAT_{\lfloor n^{\epsilon}\rfloor}(x,B(x))]$$ or $y'$ does not satisfy the assumption of $\Pr^{y'}_{x}[\cdot]_{1/n'}\ge 1/2$\big)." 
\medskip


Since $S$ is $\forall\Sigma^b_1$, by Lemma \ref{singlehard}, there is a p-time function $f_3$ and a constant $\ell$ such that \PV proves: ``$\forall 1^n>n_1, \forall m, n/2^b\le 2^{bm}\le n$, if $|h'|\ge n^{\ell}$ and a $2^{||h'||/4}$-size circuit generated by a p-time function fails to compute $h'$, then $f_3(1^n,m,h',C,x,y,y')$ outputs a circuit falsifying $$\LBtt'(f_1(1^{2^{m}}),2^{m/4},2^{m}(1/2-1/2^{m/4})),$$ or $f_3(1^n,m,h',C,x,y,y')$ outputs a circuit falsifying the assumption of $\Pr^{y'}_{x}[\cdot]_{1/n'}\ge 1/2$ or $F'$ holds," where $F'$ is the rest of the statement $S$.

Consequently, \EF proves efficiently the propositional translation of the \PV-theorem. Substituting $h'$ for $y'$, $\EF^+:=\EF+\ttable(h',2^{n/4},1/2-1/2^{n/4})$ proves the formula $F$ encoding $F'$. We now proceed as in the proof of Theorem \ref{t:core}. Assuming that $\SAT\in\Ppoly$, there is some $k$ such that for all $1^n>n_1\ge 1$ we can efficiently falsify the first disjunct of $F$. Therefore, there is a p-size circuit $B$ such that $\EF^+$ proves efficiently $\SAT_{\lfloor n^{\epsilon}\rfloor}(x,y)\rightarrow \SAT_{\lfloor n^{\epsilon}\rfloor}(x,B(x))$, which means that $\EF^+$ is p-bounded. \qed
\bigskip

\proof[Proof of Lemma \ref{l:owpwit}] 
Let $f_0(1^n,m)$ output the set of propositional formulas $A:=\{\phi_z(x)\mid z\in \{0,1\}^{cm} \}$, where $\phi_z(x)$ uses free variables $x$ together with some auxiliary variables and encodes the statement $$h(x)=h(NW_f(z)).$$ Here, $NW_f:\{0,1\}^{cm}\mapsto \{0,1\}^{\lfloor 2^{\delta m}\rfloor}$ and $c$ are given by Theorem \ref{t:nwpv}, for some $0<\gamma,\delta<1$ specified later, and $f=f_1(1^{2^m})$. The size of $\phi_z(x)$ is $\lfloor 2^{K\delta m}\rfloor$, for a constant $K$ depending only on $h$. We set $m$ so that $n/2^{K\delta}\le 2^{K\delta m}\le n$ and treat formulas $\phi_z(x)$ as formulas of size $n$. Let $A':=\{NW_{f}(z) \mid z\in \{0,1\}^{cm}\}$. 

We reason in $\SB+dWPHP(\pv)$. Suppose that an $n^k$-size circuit $C$ with $n$ inputs and $\le n$ outputs finds a satisfying assignment for all formulas in $A$. Then, there is an $n^k$-size circuit $C'$ with $u(\lfloor 2^{\delta m}\rfloor)$ inputs such that \begin{equation}\label{e:l23}2^{\lfloor 2^{\delta m}\rfloor}\times \{z\in \{0,1\}^{cm}\mid h(C'(h(NW_f(z))))\ne h(NW_f(z))\}\preceq_0 0.\end{equation}  The circuit $C'$ is obtained from $C$ by a p-time algorithm $f_{-1}$ depending on $\delta, h$ and $m$. The predicate $h(C'(h(x)))\ne h(x)$ is computable by a $2^{K'\delta m}$-size circuit $D$ with $\lfloor 2^{\delta m}\rfloor$ inputs, for a constant $K'$ depending only on $k$ and $h$. Now, we set a sufficiently small $\gamma$ and $\delta$ so that $2^{K'\delta m}\le 2^{m/4}-\lceil 2^{(\delta+\gamma)m}\rceil-1$. Therefore, by Theorem \ref{t:nwpv}, the assumption that $f$ is hard on average for $2^{m/4}$-size circuits implies that
$$2^{\lfloor 2^{\delta m}\rfloor}\times \{z<2^{cm}\mid D(NW_f(z))=1\}\succeq_{e} 2^{cm}\times \{x<2^{\lfloor 2^{\delta m}\rfloor}\mid D(x)=1\},$$ where $e=\lceil 2^{\delta m}\rceil/ 2^{m/4}$. Consequently, by (\ref{e:l23}) and Proposition \ref{lem} $ii)$, \begin{equation}\label{e:l21}0\succeq_{e} 2^{cm}\times  \{x<2^{\lfloor 2^{\delta m}\rfloor} \mid h(C'(h(x)))\ne h(x)\}.\end{equation}

We want to show that $\Pr^y_{x\in \{0,1\}^{\lfloor 2^{\delta m}\rfloor}}[h(C'(h(x)))=h(x)]_{\frac{1}{\lfloor 2^{\delta m}\rfloor}}\ge 1/2$. For the sake of contradiction, assume that this is not the case. Then, $\{x<2^{\lfloor 2^{\delta m}\rfloor} \mid h(C'(h(x))=h(x)\}\preceq_{1/\lfloor 2^{\delta m}\rfloor} 2^{\lfloor 2^{\delta m}\rfloor-1}$. By Item $1\ iii)$ of Proposition \ref{really},
$$\{x<2^{\lfloor 2^{\delta m}\rfloor} \mid h(C'(h(x))\ne h(x)\}\succeq_{2/\lfloor 2^{\delta m}\rfloor} 2^{\lfloor 2^{\delta m}\rfloor -1}.$$
By Proposition \ref{lem} $iii)$, \begin{equation}\label{e:l22}2^{cm}\times \{x<2^{\lfloor 2^{\delta m}\rfloor} \mid h(C'(h(x))\ne h(x)\}\succeq_{2/\lfloor 2^{\delta m}\rfloor} 2^{\lfloor 2^{\delta m}\rfloor-1+cm}.\end{equation} By Proposition \ref{lem} $ii)$, (\ref{e:l21}) and (\ref{e:l22}) yield $0\succeq_{2/\lfloor 2^{\delta m}\rfloor+e} 2^{\lfloor 2^{\delta m}\rfloor-1+cm}$. Hence, by Item $1\ ii)$ of Proposition \ref{really}, $2^{\lfloor 2^{\delta m}\rfloor-1+cm}<(3/\lfloor 2^{\delta m}\rfloor+e)2^{\lfloor 2^{\delta m}\rfloor-1+cm}$, which is a contradiction for 
a sufficiently big $m$. \qed

\section{Learning from the non-existence of OWFs}\label{s:dich}

\begin{theorem}[`CC $\leftarrow$ PC' from `Learning $\leftarrow$ $\not\exists$ OWF' \& hardness of {\sf E}]\label{t:dich}\ \newline
Let $k,t\ge 1$ be constants. Assume that for each sufficiently big $n$, each $2^{n/4}$-size circuit fails to compute $h'\in {\sf E}$ on $\ge 1/2-1/2^{n/4}$ of all inputs. Further, assume that there is a p-time function $h:\{0,1\}^n\mapsto\{0,1\}^{u(n)}$ such that for each constants $c,d$, there is a p-time function $f_2$ and constants $n_0$ and $0<\epsilon<1$ such that \SB proves:
\begin{itemize}
\item[] ``$\forall n$, $\forall\ cn^{c}$-size circuit $C$ with $u(m)$ inputs and $m$ outputs such that $n\le dm^d$,\\
$\big(f_2(C)$ outputs a $poly(n)$-size circuit $B$ learning $\lfloor n^{\epsilon}\rfloor^t$-size circuits with $\lfloor n^{\epsilon}\rfloor$ inputs over the uniform distribution, up to error $1/2-1/\lfloor n^{\epsilon}\rfloor$, with confidence $1/\lfloor n^{\epsilon}\rfloor$; formally, $\forall f:\{0,1\}^{\lfloor n^{\epsilon}\rfloor}\mapsto\{0,1\}$, $\forall$ $\lfloor n^{\epsilon}\rfloor^t$-size circuit $D$ computing $f$, $$\Pr^y_w[B(1^{\lfloor n^{\epsilon}\rfloor},w)\ (1/2+1/\lfloor n^{\epsilon}\rfloor)\text{-}approximates\ f]_{1/2\lfloor n^{\epsilon}\rfloor}\ge 1/\lfloor n^{\epsilon}\rfloor;$$\\ or $$\Pr^{y}_{x\in \{0,1\}^{m}}[h(C(h(x)))=h(x)]_{\frac{1}{m}}<1/2\big)."$$
\end{itemize}
Then there are constants $b$ and $a$ (depending on $k,t,h,h',c,d,f_2,n_0,\epsilon$) such that for each $n$ the existence of a function $g_n:\{0,1\}^n\mapsto \{0,1\}$ such that no circuit of size $bn^b$ computes $g_n$ on $(1/2+1/n)$ fraction of inputs and such that $\EF+\ttable(h',2^{n/4},1/2-1/2^{n/4})$ does not have $2^{an}$-size proofs of $\ttable(g_n,n^t)$ implies that 
$\SAT_n\notin\Circuit[n^{k}]$.
\end{theorem}

Note that the \SB-theorem in the assumption of Theorem \ref{t:dich} assumes $2^{\lfloor n^{\epsilon}\rfloor}\in Log$.
\medskip

\proof The assumption of Theorem \ref{t:dich} in combination with Lemma \ref{l:owpwit} implies that for any given $k,t\ge 1$, for p-time $f_1$ generating the truth-table of $h'$, there are constants $0<\epsilon<1$ and $b,n_1$ such that $\SB+dWPHP(\pv)$ proves the following statement $S$:
\medskip

\noindent  ``$\forall 1^n>n_1, \forall m, n/2^b\le 2^{bm}\le n$, if $$\LBtt'(f_1(1^{2^{m}}),2^{m/4},2^{m}(1/2-1/2^{m/4})),$$ then $f_0(1^n,m)$ outputs $A,A'\subseteq \{0,1\}^n$ such that $$\forall x\in A\ \exists y_x\in A'\ \SAT_n(x,y_x)$$ and $\big(\forall n^{k}$-size circuit $C$ with $n$ inputs and $\le n$ outputs, $$\exists x\in A, \neg\SAT_n(x,C(x))$$ or $f_2(f_{-1}(C,m))$ outputs a circuit $B$ such that 
$\forall f:\{0,1\}^{\lfloor n^{\epsilon}\rfloor}\mapsto\{0,1\}$, $\forall$ $\lfloor n^{\epsilon}\rfloor^t$-size circuit $D$ computing $f$, $$\Pr^y_w[B(1^{\lfloor n^{\epsilon}\rfloor},w)\ (1/2+1/\lfloor n^{\epsilon}\rfloor)\text{-}approximates\ f]_{1/2\lfloor n^{\epsilon}\rfloor}\ge 1/\lfloor n^{\epsilon}\rfloor,$$ or $y'$ does not satisfy the assumption of $\Pr^{y'}_{x}[\cdot]_{1/n'}\ge 1/2\big)$." 
\medskip

Since $S$ is $\forall\Sigma^b_1$, by Lemma \ref{singlehard}, there is a p-time function $f_3$ and a constant $\ell$ such that \PV proves: ``$\forall 1^n>n_1, \forall m, n/2^b\le 2^{bm}\le n$, if $|h'|\ge 2^{\ell \lfloor n^{\epsilon}\rfloor}$ and a $2^{||h'||/4}$-size circuit generated by a p-time function fails to compute $h'$, then $f_3(1^n,m,h',C,f,D,y,y')$ outputs a circuit falsifying $$\LBtt'(f_1(1^{2^{m}}),2^{m/4},2^{m}(1/2-1/2^{m/4})),$$ or $f_3(1^n,m,h',C,f,D,y,y')$ outputs a circuit falsifying the assumption of $\Pr^{y'}_{x}[\cdot]_{1/n'}\ge 1/2$ or it outputs a circuit falsifying the assumption of $\Pr^y_w[\cdot]_{1/2\lfloor n^{\epsilon}\rfloor}\ge 1/\lfloor n^{\epsilon}\rfloor$  or $F'$ holds," where $F'$ is the rest of the statement $S$.

Consequently, $\EF^+:=\EF+\ttable(h',2^{n/4},1/2-1/2^{n/4})$ proves efficiently the propositional translation of $F'$. If we now fix $n>1$ and assume that $\SAT_n\in \Circuit[n^{k''}]$, then there is some $k=O(k'')$ such that we can efficiently falsify the first disjunct of the propositional translation of $F'$ in $\EF^+$. Therefore, there is a $poly(n)$-size circuit $B$ and a $2^{Kn}$-size $\EF^+$ proof of $\mathsf{lear}^{h'}_{1/2n}(B,\Circuit[n^{t}],1/2-1/n,1/n)$, for a constant $K$ independent of $n$. Recall that this means that $\EF^+$ proves efficiently $\neg \ttable(f,n^t)\rightarrow R$, for a formula $R$.

Let $b\ge t$ be such that $B$ has size $\le bn^{b}$. We claim that for each Boolean function $g_n:\{0,1\}^n\mapsto \{0,1\}^n$ which is not $(1/2+1/n)$-approximable by any circuit of size $bn^b$, there is a $2^{an}$-size $\EF^+$-proof of $\ttable(g_n,n^{t})$, for a constant $a$ independent of $n$. This is because in order to prove $\ttable(g_n,n^t)$ in $\EF^+$, it suffices to check in $\EF^+$ that $\neg R$ holds for $f=g_n$. $\neg R$ holds for $f=g_n$ as otherwise there would be a $bn^b$-size circuit $(1/2+1/n)$-approximating $g_n$. Moreover, the fact that $\neg R$ holds for $f=g_n$ is efficiently provable in $\EF^+$ as w.l.o.g. $\neg R$, for $f=g_n$, does not contain any free variables (we can assume that the auxiliary variables are substituted by suitable constants). \qed

\section{Self-provability from random self-reducibility}\label{s:rselfred}

We show that the random self-reducibility of the discrete logarithm problem can be used to derive a conditional self-provability of the statement that the discrete logarithm problem can be solved by p-size circuits. 

For simplicity, we consider the discrete logarithm problem for $\mathbb{Z}_q^{\times}$, multiplicative groups of integers modulo a prime $q$. Let $G$ be such a cyclic group. Then there are p-time algorithms $A_1,A_2$ such that $A_1(g,h,q)=g\cdot h\in G$, for $g,h\in G$, and $A_2(g,q)=g^{-1}$, for $g\in G$. That is, $A_1$ (given $q$) defines the multiplication of two elements in $G$ and $A_2$ outputs the inverse of each $g\in G$.

The discrete logarithm problem for a cyclic group $G$ generated by $g$ is defined as follows. Given $b\in G$, we want to find $a$ such that $g^a=b$. The discrete logarithm problem is `random self-reducible': If we have a circuit $C$ which solves the problem for a $p$-fraction of all $b\in G$, we can turn it efficiently into a randomized circuit $C'$ which solves the problem on each $b\in G$ with probability $\ge p$. The circuit $C'$ interprets its random bits $r$ as $r\in [|G|]$. Then $C'$ applies $C$ on $bg^r$. Since $bg^r$ is a uniformly random element of $G$, $C$ succeeds in finding $\ell$ such that $g^{\ell}=bg^r$ with probability $\ge p$. Finally, $C'$ outputs $\ell-r$, which is the correct answer with probability $\ge p$. In other words, the following implication holds \begin{equation}\label{e:selfr1} \Pr_{b\in G}[g^{C(b)}=b]\ge p\rightarrow \forall b\in G, \Pr_{r\in [|G|]}[g^{C'(b,r)}=b]\ge p,\end{equation} where $C'$ is generated from $C$ and $g$ by a p-time function.

We want to express (\ref{e:selfr1}) by a propositional formula. To do so, we approximate probabilities by a Nisan-Wigderson generator based on a hard function $f\in\mathsf{E}$. Fix a constant $k$ and assume that $q\in (2^{n-1},2^n]$, $n\in\mathbb{N}$. Then, for each $n^k$-size circuit $C$, the predicate $g^{C(b)}=b$ with input $b$ and the predicate $g^{C'(b,r)}=b$ with input $r$ are computable respectively by circuits $D_1$ and $D_2$ with $n$ inputs and size $poly(n)$. Circuits $D_1, D_2$ reject all inputs not in $G$, so in particular $\Pr_{r\in [|G|]}[g^{C'(b,r)}=b]\le 2\Pr_{r\in \{0,1\}^n}[D_2(r)=1]$. Further, for each $\epsilon<1$, there are constants $c',c''$ and $poly(n)$-time computable generator $NW_f:\{0,1\}^{c'\lceil\log n\rceil}\mapsto \{0,1\}^n$ such that if $f:\{0,1\}^{c''\lceil\log n\rceil}\mapsto \{0,1\}$ is hard to $(1/2+1/2^{\epsilon c''\lceil\log n\rceil})$-approximate by circuits of size $2^{\epsilon c''\lceil\log n\rceil}$, then $$\left|\Pr_{z\in\{0,1\}^{c'\lceil\log n\rceil}}[D_1(NW_f(z))=1]-\Pr_{b\in \{0,1\}^n}[D_1(b)=1]\right|\le 1/n,$$ $$\left|\Pr_{z\in\{0,1\}^{c'\lceil\log n\rceil}}[D_2(NW_f(z))=1]-\Pr_{r\in \{0,1\}^n}[D_2(r)=1]\right|\le 1/n.$$

Therefore, if $f$ is hard, we have \begin{equation}\label{e:selfr2} \Pr_{z\in \{0,1\}^{c'\lceil\log n\rceil}}[D_1(NW_f(z))=1]\ge p\rightarrow \forall b\in G, \Pr_{z\in\{0,1\}^{c'\lceil\log n\rceil}}[D_2(NW_f(z))=1]\ge \frac{p}{2}-\frac{3}{2n}.\end{equation}

The advantage of (\ref{e:selfr2}) is that it can be expressed by $poly(n)$-size tautologies $\mathsf{self}_n(p,b,C)$ with
free variables for $n^k$-size circuits $C$, $n$-bit strings $b$ and $n$-bit parameters $p$ (among other extension variables). The tautologies have p-size proofs in some proof system which includes a p-size proof of the primality of $q$ and a p-size proof of the fact that $g$ is a generator of $G$. Here, we use the property that $g$ generates $G$ if and only if $g^{(q-1)/d}\not\equiv 1$ (mod $q$) for every prime $d$ dividing $q-1$. We can thus define a Cook-Reckhow propositional proof system $P_{\epsilon}$ 
 as \EF with the additional axioms which allow the system $P_{\epsilon}$ to derive any substitutional instance of $\mathsf{self}_n(p,b,C)$ and $\ttable(f,2^{\epsilon n},1/2-1/2^{\epsilon n})$ in a single step of the proof, for each sufficiently big $n$.

\begin{theorem}[Self-provability for the discrete logarithm]\label{t:selfdlog}\ \newline
Let $k$ be a constant. Assume that for some $\epsilon<1$ we have a Boolean function $f\in \mathsf{E}$ such that
for each sufficiently big $n$, $f$ is not $(1/2+1/2^{\epsilon n})$-approximable by any $2^{\epsilon n}$-size circuit. Let $P_{\epsilon}$ be the propositional proof system defined above. If there are $n^k$-size circuits solving the discrete logarithm problem for $\mathbb{Z}_q^{\times}$, where $q\in (2^{n-1},2^{n}]$, then there are p-size circuits $D$ such that $P_{\epsilon}$ has p-size proofs of tautologies encoding the statement ``$\forall b\in \mathbb{Z}_q^{\times}, g^{D(b)}=b$."
\end{theorem}

\proof
Given an $n^k$-size circuit $B$ solving the discrete logarithm problem for $\mathbb{Z}_q^{\times}$, where $q\in (2^{n-1},2^{n}]$, $P_{\epsilon}$ can derive $\mathsf{self}_n(1/2-1/n,b,B)$. Since the assumption of $\mathsf{self}_n(1/2-1/n,b,B)$ is true and since it contains essentially no free variables, it can be proven efficiently in \EF. (Essentially, in order to do so, it suffices to evaluate a \Ppoly-predicate inside \EF.) Consequently, $P_{\epsilon}$ proves efficiently $\forall b\in \mathbb{Z}_q^{\times}, g^{D(b)}=b$, for a suitable p-size circuit $D$ obtained by simulating $C'$ on all $z\in \{0,1\}^{c'\lceil\log n\rceil}$.
\qed
\bigskip


Notably, Theorem \ref{t:selfdlog} establishes a conditional equivalence between a circuit lower bound and a proof complexity lower bound (for propositional formulas which might not be tautological).

\def\zfc{
\APC-decent proof systems can be much stronger than \WF. For example, consider $\mathsf{ZFC}$ as a propositional proof system: a $\mathsf{ZFC}$-proof of propositional formula $\phi$ is a $\mathsf{ZFC}$-proof of the statement encoding that $\phi$ is a tautology. We can add the reflection of $\mathsf{ZFC}$ to \WF, i.e. we will allow \WF to derive (substitutional instances of) formulas stating that ``If $\pi$ is a $\mathsf{ZFC}$-proof of $\phi$, then $\phi$ holds." The new system is as strong as $\mathsf{ZFC}$ w.r.t. tautologies\footnote{In fact, it is equivalent to $\mathsf{ZFC}$ because $\mathsf{ZFC}$ proves efficiently its own reflection \cite{Pcon}. \cite{Pcon} also implies that $\mathsf{ZFC}$ itself is \APC-decent. } and it is easy to see that it is \APC-decent. (The reflection of the system can be proved in \APC extended with an axiom postulating the reflection for $\mathsf{ZFC}$.)}

\def\previouspaperclaim{
\begin{claim}[in \APC]\label{c:mpv} Assume that $\pi$ is a $P$-proof of $\mathsf{lear}^y_{1/2^{\ell n^{\gamma}}}(A,\Circuit[n^k],1/2-1/2^{Kn^{\gamma}},\delta)$ for a circuit $A$ and a boolean function $y$ represented by fixed bits in formula $\mathsf{lear}^y_{1/2^{\ell n^{\gamma}}}(\cdot,\cdot,\cdot,\cdot)$. Further, assume that the probability that $A$ on queries to $f$ outputs a circuit $D$ such that $\Pr[D(x)=f(x)]\ge 1/2+1/2^{Kn^{\gamma}}$ is $<\delta$, where the outermost probability is counted approximately with error $1/2^{\ell n^{\gamma}}$ using \pv-function $Sz$ and the corresponding assumptions \LBtt expressing hardness of $y$ with a suitable length $||y||=S(\cdot,\cdot,\cdot)$, i.e. using formulas $\Pr^y[\cdot]_{1/2^{\ell n^{\gamma}}}$ for the same $y$ as above - we treat $y$ as a free variable here.
 Then there is a $poly(|\pi|)$-size $P$-proof of $\ttable(f,n^k)$ or $y$ does not satisfy the assumptions of $\Pr^y[\cdot]_{1/2^{\ell n^{\gamma}}}$.
\end{claim}

To see that the claim holds, we reason in \APC as follows. Assume $\pi$ is a $P$-proof of $\mathsf{lear}^y_{1/2^{\ell n^{\gamma}}}(A,\Circuit[n^k],1/2-1/2^{Kn^{\gamma}},\delta)$ but $A$ on queries to $f$ outputs a circuit $(1/2+1/2^{Kn^{\gamma}})$-approximating $f$ with probability $<\delta$. Then, either $y$ does not satisfy the assumptions of $\Pr^y[\cdot]_{1/2^{\ell n^{\gamma}}}$ or there is a trivial $2^{O(n)}$-size $P$-proof of $\neg\ttable(f,n^k)\rightarrow\neg R(b)$, for predicate $R$ from the definition of $\mathsf{lear}^y_{1/2^{\ell n^{\gamma}}}(A,\Circuit[n^k],1/2-1/2^{Kn^{\gamma}},\delta)$ and a complete assignment $b$. The $P$-proof is obtained by evaluating function $Sz$ which counts the confidence of $A$ - note that functions $f,y$ and algorithm $A$ are represented inside $P$ by fixed bits so the $P$-proof just evaluates a $2^{O(n)}$-size circuit on some input, which is possible by Lemma \ref{l:decent}, Item 2. (We use here also the fact that \APC knows that the probability statement expressed by function $Sz$ translates to $\neg R$ in the negation normal form.) 
 The formula $\neg\ttable(f,n^k)\rightarrow\neg R(b)$ is obtained from $\neg R(b)$ by an instantiation of a single Frege rule, which is available by Lemma \ref{l:decent}, Item 1. 
Applying again Lemma \ref{l:decent}, Item 1, from a $P$-proof of $\mathsf{lear}^y_{1/2^{\ell n^{\gamma}}}(A,\Circuit[n^k],1/2-1/2^{Kn^{\gamma}},\delta)$ and a $P$-proof of $\neg\ttable(f,n^k)\rightarrow\neg R(b)$, we construct a $poly(|\pi|)$-size $P$-proof of $\ttable(f,n^k)$. This proves the claim.
}

\def\previouspaper{
\begin{lemma}[from \cite{BFKL}]\label{l:gen}
There is a randomized p-time function $L$ such that for every $n^c$-size circuit $C$, 
if an $s$-size circuit $D$ satisfies $$\Pr[D(x)=1]-\Pr[D(G_C(x))=1]\ge 1/s,$$ then the circuit $C$ is learnable by $L(D)$ over the uniform distribution with random examples, confidence $1/2m^2s$, up to error $1/2-1/2ms$.
\end{lemma}

\proof Given $D$, $L(D)$ chooses a random $i\in [m]$, random bits $r_{i},\dots,r_m$, random $n$-bit strings $x_1,\dots,x_n$ except $x_i$ and queries the bits $C(x_1),\dots, C(x_{i-1})$. For $x_i\in\{0,1\}^n$, let $p_i:=D(x_1,C(x_1),\dots,x_{i-1},C(x_{i-1}),x_i,r_i,\dots,x_m,r_m)$. Then $L(D)$ on $x_i$ predicts the value $C(x_i)$ by outputting $\neg r_i$ if $p_i=1$ and $r_i$ otherwise. By triangle inequality, random $i\in [m]$ satisfies $$\Pr[p_{i}=1]-\Pr[p_{i+1}=1]\ge 1/ms$$ with probability $1/m$. Since the probability over $r_i\dots,r_m,x_1,\dots,x_m$ that $L(D)$ predicts $C(x_i)$ correctly is $$\frac{1}{2}\Pr[p_i=1\mid r_i\ne C(x_i)]+\frac{1}{2}(1-\Pr[p_i=1\mid r_i= C(x_i)]),$$ and $\Pr[p_i=1]=\frac{1}{2} \Pr[p_i=1\mid r_i=C(x_i)]+\frac{1}{2}\Pr[p_i=1\mid r_i\ne C(x_i)],$ it follows that $$\Pr_{x_i}[L(D)(x_{i})=C(x_i)]\ge 1/2+1/2ms$$ with probability $1/2m^2s$ over the internal randomness of $L(D)$. \qed
\medskip

\begin{question}[Dichotomy]\label{q:dichotomy} Assume that for each $\epsilon<1$ there is no pseudorandom generator $g:\{0,1\}^n\mapsto \{0,1\}^{n+1}$ computable in \Ppoly and safe against circuits of size $2^{n^{\epsilon}}$ for infinitely many $n$. Does it follow that p-size circuits are learnable by circuits of size $2^{O(n^\delta)}$, for some $\delta<1$, with confidence $1/n$, up to error $1/2-1/2^{O(n^\delta)}$?
\end{question}

\subsection{Worst-case learning from strong lower bound methods} 

\begin{definition} The circuit size problem $\GCSP[s,k]$ is the problem to decide whether for a given list of $k$ samples $(y_i,b_i)$, $y_i\in\{0,1\}^n, b_i\in\{0,1\}$, there exists a circuit $C$ of size $s$ computing the partial function defined by samples $(y_i,b_i)$, i.e. $C(y_i)=b_i$ for the given $k$ samples $(y_i,b_i)$. The parameterized minimum circuit size problem $\MCSP[s]$ stands for $\GCSP[s,2^n]$ where the list of $2^n$ samples defines the whole truth-table of a Boolean function.
\end{definition}

\begin{theorem}[Learning from succinct natural proofs]\label{t:gcsplear}
Assume $\GCSP[n^c,n^{d}]\in\Ppoly$ for constants $d>c+1$. Then, $\Circuit[n^c]$ is learnable by \Ppoly over the uniform distribution with random examples, confidence $1/poly(n)$, up to error $1/2-1/poly(n)$.
\end{theorem}

\proof As the number of partial Boolean functions on a given set of $m$ inputs is $2^m$ and the number of $n^c$-size circuits is bouded by $2^{n^{c+1}}$, 
$\GCSP[n^c,n^{d}]\in\Ppoly$ implies that for $m=n^{d}$ there are p-size circuits $D$ such that for each $n^c$-size circuit $C$, $$\Pr[D(x)=1]-\Pr[D(G_C(x))=1]\ge 1/2.$$ Now, it suffices to apply Lemma \ref{l:gen}. \qed 

\subsection{Learning from breaking pseudorandom function families}\label{s:prfs}

Oliveira and Santhanam \cite{OS} showed that the assumption of the existence of natural proofs from Theorem \ref{t:cikk} can be further weakened to the existence of a distinguisher breaking non-uniform pseudorandom function families. Their result follows from a combination of Theorem \ref{t:cikk} and the Min-Max Theorem. Using their strategy but combining the Min-Max Theorem with Theorem \ref{t:gcsplear}, learning algorithms with random examples can be obtained from distinguishers breaking succinct non-uniform pseudorandom function families
\bigskip

A {\em two-player zero-sum game} is specified by an $r\times c$ matrix $M$ and is played as follows. MIN, the row player, chooses a probability distribution $p$ over the rows. MAX, the column player, chooses a probability distribution $q$ over the columns. A row $i$ and a column $j$ are drawn randomly from $p$ and $q$, and MIN pays $M_{i,j}$ to MAX. MIN plays to minimize the expected payment, MAX plays to maximize it. The rows and columns are called the {\em pure strategies} available to MIN and MAX, respectively, while the possible choices of $p$ and $q$ are called {\em mixed strategies}. The Min-Max theorem states that playing first and revealing one's mixed strategy is not a disadvantage: 
$$min_{p} max_j \sum_i p(i)M_{i,j}=max_q min_i \sum_j q(j) M_{i,j}.$$ Note that the second player need not play a mixed strategy - once the first player's strategy is fixed, the expected payoff is optimized for the second player by playing some pure strategy. The expected payoff when both players play optimally is called the {\em value} of the game. We denote it $v(M)$.

A mixed strategy is {\em $k$-uniform} if it chooses uniformly from a multiset of $k$ pure strategies. Let $M_{min}=min_{i,j} M_{i,j}$ and $M_{max}=max_{i,j} M_{i,j}$. Newman \cite{Nm}, Alth\"ofer \cite{Alt} and Lipton-Young \cite{LY} showed that each player has a near-optimal $k$-uniform strategy for $k$ proportional to the logarithm of the number of pure strategies available to the opponent.

\begin{theorem}[\cite{Nm, Alt, LY}]\label{thm:minmax} For each $\epsilon>0$ and $k\ge \ln(c)/2\epsilon^2$, $$min_{p\in P_k} max_j \sum_{i} p(i) M_{i,j}\le v(M)+\epsilon(M_{max}-M_{min}),$$ where $P_k$ denotes the $k$-uniform strategies for MIN. The symmetric result holds for MAX.
\end{theorem}

\begin{definition}[Succinct non-uniform PRF] An $(m,m')$-succinct non-uniform pseudorandom function family from circuit class $\mathcal{C}$ safe against circuits of size $s$ is a set $S$ of partial truth-tables $\langle (x_1,b_1),\dots, (x_m,b_m)\rangle$ where each $x_i$ is an $n$-bit string and $b_i\in\{0,1\}$ such that each partial truth-table from $S$ is computable by one of $m'$ circuits from $\mathcal{C}$ and for every circuit $D$ of size $s$, $$\Pr_{x}[D(x)=1]-\Pr_{x\in S}[D(x)=1]<1/s$$ where the first probability is taken over $x\in\{0,1\}^{m(n+1)}$ chosen uniformly at random and the second probability over partial truth-tables chosen uniformly at random from $S$.
\end{definition}

\begin{theorem}[Learning or succinct non-uniform PRF]\label{t:towardscore}
Let $c\ge 1$ and $s>n,m\ge 1$. There is an $(m,8s^4)$-succinct non-uniform PRF in $\Circuit[n^c]$ safe against $\Circuit[s]$ or there are circuits of size $poly(s)$ learning $\Circuit[n^c]$ over the uniform distribution with random examples, confidence $1/poly(s)$, up to error $1/2-1/poly(s)$.
\end{theorem}

\proof Consider a two-player zero-sum game specified by a matrix $M$ with rows indexed by $n^c$-size circuits with $n$ inputs and columns indexed by $s$-size circuits with $m(n+1)$ inputs. Define the entry $M_{C,D}$ of $M$ corresponding to a row circuit $C$ and a column circuit $D$ as $$M_{C,D}:=|\Pr_x[D(x)=1]-\Pr_x[D(G_C(x))=1]|$$ for the generator $G_C$ from the proof of Lemma \ref{l:gen}. Hence $M_{max}-M_{min}\le 1$. 

If $v(M)\ge 1/4s$, then by Theorem \ref{thm:minmax} (with $\epsilon=1/8s$), there exist a multiset of $k\le 32n^{c+1}s^2$ $s$-size circuits $D^1,\dots, D^{k}$ such that for every $n^c$-size circuit $C$, a random $D$ from $D^1,\dots, D^{k}$ satisfies $$\text{E}[|\Pr[D(x)=1]-\Pr[D(G_C(x))=1]|]\ge 1/8s.$$ 

By Lemma \ref{l:gen}, for every $n^c$-size circuit $C$, one of the circuits $D^1,\dots, D^k$ (or their negations) can be used to learn $C$ with confidence $1/poly(s)$, up to error $1/2-1/poly(s)$. 
A $poly(s)$-size circuit using a random $D^i$ from $D^1,\dots, D^k$ or its negation thus learns $\Circuit[n^c]$ with random examples, confidence $1/poly(s)$, up to error $1/2-1/poly(s)$.
\medskip

If $v(M)<1/4s$, then by Theorem \ref{thm:minmax} (with $\epsilon=1/4s$), there exists a multiset of $k\le 8s^4$ $n^c$-size circuits $C^1,\dots, C^k$ such that for every $s$-size circuit $D$, a random $C$ from $C^1,\dots, C^k$ satisfies $$\text{E}[|\Pr[D(x)=1]-\Pr[D(G_C(x))=1]|]\le 1/2s.$$ Since $\text{E}[|\Pr[D(x)=1]-\Pr[D(G_C(x))=1]|]\ge|\Pr[D(x)=1]-\text{E}[\Pr[D(G_C(x))=1]]|$ a generator $$G:\{0,1\}^{mn+\lceil\log k\rceil}\mapsto \{0,1\}^{mn+m}$$ which takes as input a string of length $mn+\lceil\log k\rceil$ encoding (an index of) a circuit $C$ from $C^1,\dots, C^k$ together with $m$ $n$-bit strings $x_1,\dots,x_m$ and outputs $x_1,C(x_1),\dots, x_m, C(x_m)$ is safe against circuits of size $s$. The range of $G$ defines an $(m,8s^4)$-succinct non-uniform PRF in $\Circuit[n^c]$ safe against $\Circuit[s]$.
\qed
}

\section*{Acknowledgement}

J\'an Pich received support from the Royal Society University Research Fellowship URF$\backslash$R1$\backslash$211106 ``Proof complexity and circuit complexity: a unified approach." Rahul Santhanam was partially funded by the EPSRC New Horizons grant EP$\backslash$V048201$\backslash$1: ``Structure versus Randomness in Algorithms and Computation." Part of this work was done while the authors were visiting the Simons Institute for the Theory of Computing. For the purpose of Open Access, the authors have applied a CC BY public copyright lincence to any Author Accepted Manuscript version arising from this submission.


\begin{thebibliography}{00}


\bibitem{Alt} Alth\"ofer I.; {\it On sparse approximations to randomized strategies and convex combinations}; Linear Algebra and its Applications, 199(1):339-355, 1994.




\bibitem{Abb} Atserias A.; {\it Distinguishing SAT from polynomial-size circuits, through black-box queries}; Computational Complexity Conference (CCC), 2006.

\bibitem{BP} Beame P., Pitassi T.; {\it Propositional Proof Complexity: Past, Present and Future}; Current Trends in Theoretical Computer Science, Entering the 21st Century, 42-70, 2001.

\bibitem{BCKRS} Binnendyk E., Carmosino M., Kolokolova A., Ramyaa R., Sabin M.; {\it Learning with distributional inverters}; Algorithmic Learning Theory (ALT), 2022.


\bibitem{Hardins} Bogdanov A., Talwar K., Wan A.; {\it Hard instances for satisfiability and quasi-one-way functions}; ICS, 2010.



\bibitem{Bba} Buss S.; {\it Bounded arithmetic}; Bibliopolis, 1986.










\bibitem{Cptime} Cobham A.; {\it The intrinsic computational difficulty of functions}; Proceedings of the 2nd International Congress of Logic, Methodology and Philosophy of Science, North Holland, pp. 24-30, 1965.

\bibitem{Cpv} Cook S.A.; {\it Feasibly constructive proofs and the propositional calculus}; Symposium on Theory of Computing (STOC), 1975.


\bibitem{CTcol} Cook S.A., Thapen N.; {\it The strength of replacement in weak arithmetic}; ACM Transactions on Computational Logic, 7(4):749-764, 2006. 

\bibitem{CK} Cook S.A., Kraj\'{i}\v{c}ek J.; {\it Consequences of the Provability of NP$\subseteq$P/poly}; Journal of Symbolic Logic, 72:1353-1357, 2007.

\bibitem{Ls} de Rezende S., G\"o\"os M., Robere R.; {\it Proofs, Circuits and Communication}; SIGACT News Complexity Theory Column, 2022.




\bibitem{GP} Grochow J., Pitassi T.; {\it Circuit Complexity, Proof Complexity, and Polynomial Identity Testing: The Ideal Proof System}; Journal of the ACM, 65(6), 37:1-59, 2018.

\bibitem{GST} Gutfreund D., Shaltiel R., Ta-Shma A.; {\it If NP languages are hard in the worst-case then it is easy to find their hard instances}; Computational Complexity, 16(4), 412-441, 2007.

\bibitem{Hb} Hirahara S.; {\it Non-black-box worst-case to average-case reductions within \NP}; Foundations of Computer Science (FOCS), 2018.


\bibitem{Iw} Impagliazzo R.; {\it A personal view of average-case complexity}; Structure in Complexity Theory (SCT), 1995.



\bibitem{Jwphp} Je\v{r}\'abek E.; {\it Dual weak pigeonhole principle, Boolean complexity and derandomization}; Annals of Pure and Applied Logic, 129:1-37, 2004.

\bibitem{Jphd} Je\v{r}\'abek E.; {\it Weak pigeonhole principle and randomized computation}; Ph.D. thesis, Charles University in Prague, 2005.

\bibitem{Japx} Je\v{r}\'abek E.; {\it Approximate counting in bounded arithmetic}; Journal of Symbolic Logic, 72:959-993, 2007.

\bibitem{Klb} Kannan R.; {\it Circuit-size lower bounds and non-reducibility to sparse sets}, Information and Control, 55:40-56, 1982. 

\bibitem{Kba} Kraj\'{i}\v{c}ek J.; {\it Bounded arithmetic, propositional logic, and complexity theory}; Cambridge University Press, 1995.


\bibitem{Kwphp} Kraj\'{i}\v{c}ek J.; {\it On the weak pigeonhole principle}; Fundamenta Mathematicae, 170(1-3):123-140, 2001.





\bibitem{Kpc} Kraj\'{i}\v{c}ek J.; {\it Proof complexity}; Cambridge University Press, 2019.


\bibitem{KPT} Kraj\'{i}\v{c}ek J., Pudl\'ak P., Takeuti G.; {\it Bounded arithmetic and the polynomial hierarchy}; Annals of Pure and Applied Logic, 52:143-153, 1991.





\bibitem{LY} Lipton R.J., Young N.E.; {\it Simple strategies for large zero-sum games with applications to complexity theory}; Symposium on Theory of Computing (STOC), 1994.

\bibitem{LP} Liu Y., Pass R.; {\it On one-way functions from \NP-complete problems}; Computational Complexity Conference (CCC), 2022.



\bibitem{MP} M\"uller M., Pich J.; 
{\it Feasibly constructive proofs of succinct weak circuit lower bounds}; Annals of Pure and Applied Logic, 2019.  

\bibitem{Nm} Newman I.; {\it Private vs common random bits in communication complexity}; Information Processing Letters, 39:67-71, 1991.








\bibitem{Pclba} Pich J.; {\it Circuit lower bounds in bounded arithmetics}; Annals of Pure and Applied Logic, 166(1):29-45, 2015.



\bibitem{PSlearaut} Pich J., Santhanam R.; {\it Learning algorithms versus automatability of Frege systems}; arXiv, 2021.




\bibitem{Rup} Razborov A.A.; {\it Bounded arithmetic and lower bounds in boolean complexity}; Feasible Mathematics II, 344-386, 1995. 


\bibitem{Rkdnf} Razborov A.A.; {\it Pseudorandom generators hard for $k$-DNF Resolution and Polynomial Calculus}; Annals of Mathematics, 181(2):415-472, 2015. 



\bibitem{S19} Santhanam R.; {\it Pseudorandomness and the Minimum Circuit Size Problem}; Innovations in Theoretical Computer Science (ITCS), 2020.




\bibitem{Wacc} Williams R.; {\it Non-uniform ACC circuit lower bounds}; Computational Complexity Conference (CCC), 2011.

\end{thebibliography}
\end{document}